\documentclass[11pt]{article}
\usepackage{mathrsfs}
\usepackage{amsthm}
\usepackage{amssymb}
\usepackage{amsmath}
\usepackage{graphicx}
\usepackage{color}
\usepackage{amsfonts}
\usepackage{float}
\usepackage{cite}
\usepackage [latin1]{inputenc}
\usepackage[text={152mm,228mm},left=38mm,vmarginratio=1:1]{geometry}
\newtheorem{theorem}{Theorem}[section]

\newtheorem{lemma}[theorem]{Lemma}

\numberwithin{equation}{section}
\normalsize

\begin{document}
\title{\textbf{Mean field limits of a class of conservative systems with position-dependent transition rates}}

\author{Xiaofeng Xue \thanks{\textbf{E-mail}: xfxue@bjtu.edu.cn \textbf{Address}: School of Mathematics and Statistics, Beijing Jiaotong University, Beijing 100044, China.}\\ Beijing Jiaotong University}

\date{}
\maketitle

\noindent {\bf Abstract:} In this paper, we are concerned with a class of conservative systems including asymmetric exclusion processes and zero-range processes as examples, where some particles are initially placed on $N$ positions. A particle jumps from a position to another at a rate depending on coordinates of these two positions and numbers of particles on these two positions. We show that the hydrodynamic limit of our model is driven by a nonlinear function-valued ordinary differential equation which is consistent with a mean field analysis. Furthermore, in the case where numbers of particles on all positions are bounded by $\mathcal{K}<+\infty$, we show that the fluctuation of our model is driven by a generalized Ornstein-Uhlenbeck process. A crucial step in proofs of our main results is to show that numbers of particles on different positions are approximately independent by utilizing a graphical method.

\quad


\section{Introduction}\label{section one}

In this paper, we are concerned with a class of conservative systems on complete graphs with position-dependent transition rates, which include zero range processes and generalized exclusion processes as special cases. For later use, we first introduce some nations and definitions. We denote by $\overline{\mathbb{Z}}_+$ the set $\{0,1,2,\ldots\}\bigcup \{+\infty\}$. For any $\mathcal{K}\in \overline{\mathbb{Z}}_+$, we denote by $\mathbb{Z}_\mathcal{K}$ the set $\{0,1,2,\ldots, \mathcal{K}\}$ when $\mathcal{K}<+\infty$ and the set $\{0,1,2,\ldots\}$ when $\mathcal{K}=+\infty$. We denote by $\mathbb{T}$ the one-dimensional torus $[0, 1)$. For a given integer $N\geq 1$, our conservative system $\{\eta_t\}_{t\geq 0}$ on the complete graph of $N$ vertices is a continuous-time Markov process with state space $(\mathbb{Z}_\mathcal{K})^N$ for some $\mathcal{K}\in \overline{\mathbb{Z}}_+$, i.e., for any $t\geq 0$, $\eta_t$ is a configuration which represents some particles being placed on $N$ positions and $\eta_t(j)$ is the number of particles on the $j$th position at the moment $t$ for $1\leq j\leq N$. Our process $\{\eta_t\}_{t\geq 0}$ evolves as follows. For $1\leq i\neq j \leq N$ and $k,l\in \mathbb{Z}_\mathcal{K}$, if there are $k$ particles on the $i$th position and $l$ particles on the $j$th position, then a particle jumps from the $i$th position to the $j$th position at rate $\frac{1}{N}\phi_{k,l}\left(\frac{i}{N}, \frac{j}{N}\right)$, where $\phi_{k,l}\in C^\infty(\mathbb{T}^2)$ are nonnegative for all $k,l\in \mathbb{Z}_{\mathcal{K}}$. Since the number of particles is nonnegative, we further define $\phi_{0,l}\equiv 0$ for all $l\in \mathbb{Z}_{\mathcal{K}}$. When $\mathcal{K}<+\infty$, each position has at most $\mathcal{K}$ particles, hence we further define $\phi_{k, \mathcal{K}}\equiv0$ for all $k\in \mathbb{Z}_{\mathcal{K}}$ when $\mathcal{K}<+\infty$.

The process $\{\eta_t\}_{t\geq 0}$ can be defined equivalently according to its generator. For any $\eta\in (\mathbb{Z}_\mathcal{K})^N$ and $1\leq i\neq j\leq N$, we define $\eta^{(i,j)}\in (\mathbb{Z}_{\mathcal{K}})^N$ as
\[
\eta^{(i,j)}(x)=
\begin{cases}
\eta(x) & \text{~if~}x\neq i,j,\\
\eta(i)-1 & \text{~if~}x=i,\\
\eta(j)+1 & \text{~if~}x=j.
\end{cases}
\]
According the evolution of $\{\eta_t\}_{t\geq 0}$ given above, the generator $\mathcal{L}$ of $\{\eta_t\}_{t\geq 0}$ is given by
\begin{equation}\label{equ 1.1 generator}
\mathcal{L}f(\eta)=\frac{1}{N}\sum_{i=1}^N\sum_{j\neq i}\sum_{k\in \mathbb{Z}_\mathcal{K}}\sum_{l\in \mathbb{Z}_{\mathcal{K}}}
\phi_{k,l}\left(i/N, j/N\right)1_{\{\eta(i)=k, \eta(j)=l\}}\left(f(\eta^{(i,j)})-f(\eta)\right)
\end{equation}
for any $\eta\in (\mathbb{Z}_\mathcal{K})^N$ and $f$ from $(\mathbb{Z}_\mathcal{K})^N$ to $\mathbb{R}$, where $1_A$ is the indicator function of the event $A$. From now on, to emphasize the dependence of $N$, we replace notations $\eta_t$ and $\mathcal{L}$ by $\eta_t^N$ and $\mathcal{L}_N$ respectively. We denote by $\mathbb{P}$ the probability measure of our processes and by $\mathbb{E}$ the expectation with respect to $\mathbb{P}$.

Our processes include following important special cases.

\quad

\textbf{Example 1} (\emph{Generalized exclusion process}). When $\mathcal{K}<+\infty$, our model reduces to the generalized exclusion process (see Section 2.4 of \cite{kipnis+landim99}), where each position has at most $\mathcal{K}$ particles and a particle jumps at a rate depending on numbers of particles on the starting and the ending positions. Especially, when $\mathcal{K}=1$, our model reduces to the exclusion process (see Chapter 8 of \cite{Lig1985} or Part {\rm \uppercase\expandafter{\romannumeral3}} of \cite{Lig1999}), where each position has at most $1$ particle and all particles perform random walks where each jump to an already occupied position is suppressed.

\quad

\textbf{Example 2} (\emph{Zero range process}). When $\mathcal{K}=+\infty$ and $\phi_{k,l}$ is independent of $l$ for all $k,l$ and hence $\phi_{k,l}$ can be rewritten as $\phi_k$, our model reduces to the zero-range process (see Section 2.3 of \cite{kipnis+landim99}), where a particle jumps at a rate in dependent of the number of particles on the ending position. Especially, when $\phi_k=k\phi$ for some $\phi\in C^\infty(\mathbb{T}^2)$, our model reduces to the generalized $N$-urn Ehrenfest model (see \cite{Xue2023}), where each particle independently jumps at a rate only depending on coordinates of the starting and ending positions.

\quad

\textbf{Example 3} (\emph{Misanthrope process}). For any given $(u,v)\in \mathbb{T}^2$, if $\phi_{k,l}(u,v)$ is increasing in $k$ and decreasing in $l$, then our model is also called a misanthrope process (see \cite{Cocozza1985}), where particles prefer to jump to positions with less other particles. Note that this example has overlaps with above two. The exclusion process and the Ehrenfest model are both misanthrope processes.

\quad

In this paper, we are concerned with the mean field limit of the empirical density field of $\{\eta_t^N\}_{t\geq 0}$ as $N\rightarrow+\infty$, which is also called as the hydrodynamic limit, i.e., we will give the law of large numbers of $\frac{1}{N}\sum_{j=1}^N1_{\{\eta_t^N(j)=k\}}f(j/N)$ for all $k\in \mathbb{Z}_{\mathcal{K}}$ and $f\in C(\mathbb{T})$ as $N\rightarrow+\infty$. We will show that the hydrodynamic limit is driven by a $C(\mathbb{T})$-valued nonlinear ordinary differential equation (ODE), which is consistent with a heuristic mean-field analysis. The central limit theorem from the above hydrodynamic limit, which is also called as the fluctuation, is discussed for the case where $\mathcal{K}<+\infty$. We will show that the fluctuation is driven by a generalized Ornstein-Uhlenbeck (O-U) process. For mathematical details, see Section \ref{section two}.

Hydrodynamics of several types of interacting particle systems on complete graphs with position-dependent transition rates, which are also called as $N$-urn interacting particle systems, are discussed in previous literature. Reference \cite{Xue2023} investigates hydrodynamics of the $N$-urn Ehrenfest model. It is shown in \cite{Xue2023} that the hydrodynamic limit of the $N$-urn Ehrenfest model is driven by a linear $C(\mathbb{T})$-valued ODE and the corresponding fluctuation is driven by a generalized O-U process. Reference \cite{Xue2023b} introduces a class of $N$-urn linear systems, which include voter models, symmetric exclusion processes and binary contact path processes as special cases. Similar with results given \cite{Xue2023}, it is proved in \cite{Xue2023b} that the hydrodynamic limit of the $N$-urn linear system is driven by a linear $C(\mathbb{T})$-valued ODE and the fluctuation is driven by a O-U process.

For the proof of the hydrodynamic limit of a $N$-urn interacting particle system, a crucial step is to show that states of different positions are approximately independent as $N\rightarrow+\infty$, where different models require different strategies. For the Ehrenfest model, numbers of particles on different positions are independent at any moment $t>0$ when at $t=0$ the number of particles on each position independently follows a Poisson distribution. For linear systems, the above approximated independence are proved by directly estimating both density and correlation fields of the process via a linear-system version of Chapman-Kolmogorov equation introduced in \cite{Lig1985}. In this paper, we show that numbers of particles on different positions are approximately independent according to a graphical method. For mathematical details, see Section \ref{section three}.

\section{Main results}\label{section two}

In this section, we give our main results. For later use, we first introduce some notations and assumptions. For any $\mathcal{K}\in \overline{\mathbb{Z}}_+$, if $\mathcal{K}<+\infty$, we denote by $\mathcal{A}_\mathcal{K}$ the set of functions $\psi\in C^\infty(\mathbb{T})$ such that $0<\psi(u)<\mathcal{K}$ for all $u\in \mathbb{T}$. If $\mathcal{K}=+\infty$, we denote by $\mathcal{A}_{+\infty}$ the set of functions $\psi\in C^\infty(\mathbb{T})$ such that $\psi(u)>0$ for all $u\in \mathbb{T}$. For any $N\geq 1$, $\mathcal{K}\in \overline{\mathbb{Z}}_+$ and $\psi\in \mathcal{A}_\mathcal{K}$, if $\mathcal{K}<+\infty$, then we denote by $\nu^N_{\psi,\mathcal{K}}$ the probability measure on $(\mathbb{Z}_\mathcal{K})^N$ under which $\{\eta(i)\}_{1\leq i\leq N}$ are independent and $\eta(i)$ follows the binomial distribution $B(\mathcal{K},\psi(i/N)/\mathcal{K})$ for all $1\leq i\leq N$. If $\mathcal{K}=+\infty$, then we denote by $\nu^N_{\psi,+\infty}$ the probability measure on $(\mathbb{Z}_{+\infty})^N$ under which $\{\eta(i)\}_{1\leq i\leq N}$ are independent and $\eta(i)$ follows the Poisson distribution with parameter $\psi(i/N)$ for all $1\leq i\leq N$. For any $t\geq 0$ and $k\in \mathbb{Z}_\mathcal{K}$, we denote by $\mu_{t,k}^N(du)$ the random measure
\[
\frac{1}{N}\sum_{i=1}^N1_{\{\eta_t^N(i)=k\}}\delta_{i/N}(du)
\]
on $\mathbb{T}$, where $\delta_a(du)$ is the Dirac measure concentrated at $a$. We can also consider $\mu_{t,k}(du)$ as a random element in $(C(\mathbb{T}))^\prime$ in the sense that
\[
\mu_{t,k}^N(f)=\frac{1}{N}\sum_{i=1}^N1_{\{\eta_t^N(i)=k\}}f(i/N)
\]
for any $f\in C(\mathbb{T})$. Throughout this paper, we adopt the following assumption of the initial distribution.

\textbf{Assumption (A)}. There exists $\psi\in \mathcal{A}_\mathcal{K}$ such that the distribution of $\eta^N_0$ is $\nu^N_{\psi, \mathcal{K}}$ for all $N\geq 1$, where $\mathcal{K}\in \overline{\mathbb{Z}}_+$ is the parameter of our model defined as in Equation \eqref{equ 1.1 generator}.

For $\mathcal{K}<+\infty$, we are concerned with the following $\left(C(\mathbb{T})\right)^{\mathcal{K}+1}$-valued ordinary differential equation. For any $u\in \mathbb{T}$ and $t\geq 0$,

\begin{equation}\label{equ 2.1 mean field ODE finite types}
\begin{cases}
&\frac{d}{dt}\rho_{t,k}^\mathcal{K}(u)=-\rho^\mathcal{K}_{t,k}(u)\sum_{l=0}^\mathcal{K}\int_{\mathbb{T}}\phi_{k,l}(u,v)\rho_{t,l}^\mathcal{K}(v)dv
-\rho_{t,k}^{\mathcal{K}}(u)\sum_{l=1}^\mathcal{K}\int_{\mathbb{T}}\phi_{l,k}(v,u)\rho_{t,l}^{\mathcal{K}}(v)dv\\
&\text{\quad\quad}+\rho^{\mathcal{K}}_{t,k-1}(u)\sum_{l=1}^\mathcal{K}\int_{\mathbb{T}}\phi_{l,k-1}(v,u)\rho^{\mathcal{K}}_{t,l}(v)dv
+\rho_{t,k+1}^{\mathcal{K}}(u)\sum_{l=0}^\mathcal{K}\int_{\mathbb{T}}\phi_{k+1,l}(u,v)\rho_{t,l}^\mathcal{K}(v)dv\\
&\text{\quad\quad if~}1\leq k\leq \mathcal{K}-1,\\
&\frac{d}{dt}\rho_{t,0}^\mathcal{K}(u)=\rho_{t,1}^{\mathcal{K}}(u)\sum_{l=0}^\mathcal{K}\int_{\mathbb{T}}\phi_{1,l}(u,v)\rho_{t,l}^\mathcal{K}(v)dv
-\rho_{t,0}^{\mathcal{K}}(u)\sum_{l=1}^\mathcal{K}\int_{\mathbb{T}}\phi_{l,0}(v,u)\rho_{t,l}^{\mathcal{K}}(v)dv,\\
&\frac{d}{dt}\rho_{t,\mathcal{K}}^\mathcal{K}(u)
=\rho^{\mathcal{K}}_{t,\mathcal{K}-1}(u)\sum_{l=1}^\mathcal{K}\int_{\mathbb{T}}\phi_{l,\mathcal{K}-1}(v,u)\rho^{\mathcal{K}}_{t,l}(v)dv
-\rho^\mathcal{K}_{t,\mathcal{K}}(u)\sum_{l=0}^\mathcal{K}\int_{\mathbb{T}}\phi_{\mathcal{K},l}(u,v)\rho_{t,l}^\mathcal{K}(v)dv, \\
&\rho_{0,k}^\mathcal{K}(u)={\mathcal{K} \choose k}\left(\psi(u)/\mathcal{K}\right)^k\left(1-\frac{\psi(u)}{\mathcal{K}}\right)^{\mathcal{K}-k}
\text{~for~}0\leq k\leq \mathcal{K}.
\end{cases}
\end{equation}

Now we state our first main result, which gives the hydrodynamic limit of our model in the case where $\mathcal{K}<+\infty$.

\begin{theorem}\label{theorem 2.1 mean field limit finite}
If $\mathcal{K}<+\infty$, then the solution $\{\rho^{\mathcal{K}}_{t,k}\}_{0\leq k\leq \mathcal{K}}$ to Equation \eqref{equ 2.1 mean field ODE finite types} exists for $t\in [0, +\infty)$ and is unique. Furthermore, if $\mathcal{K}<+\infty$, then, under Assumption (A),
\begin{equation}\label{equ 2.2 convergence in probability finite}
\lim_{N\rightarrow+\infty}\mu_{t,k}^N(f)=\int_{\mathbb{T}}\rho_{t,k}^\mathcal{K}(u)f(u)du
\end{equation}
in $L^2$ for any $t\geq 0$, $0\leq k\leq \mathcal{K}$ and $f\in C(\mathbb{T})$.
\end{theorem}

To give an analogue result of Theorem \ref{theorem 2.1 mean field limit finite} in the case where $\mathcal{K}=+\infty$, we need a further assumption.

\textbf{Assumption (B)}. When $\mathcal{K}=+\infty$, then there exists $C_1<+\infty$ such that
\[
\phi_{k,l}(u,v)\leq C_1k
\]
for any $0\leq k,l<+\infty$ and $u, v\in \mathbb{T}$.

For $\mathcal{K}=+\infty$, we are concerned with the following $\left(C(\mathbb{T})\right)^\infty$-valued ordinary differential equation. For any $u\in \mathbb{T}$ and $t\geq 0$,

\begin{equation}\label{equ 2.3 mean field ODE infinite types}
\begin{cases}
&\frac{d}{dt}\rho_{t,k}^\infty(u)=-\rho^\infty_{t,k}(u)\sum_{l=0}^\infty\int_{\mathbb{T}}\phi_{k,l}(u,v)\rho_{t,l}^\infty(v)dv
-\rho_{t,k}^\infty(u)\sum_{l=1}^\infty\int_{\mathbb{T}}\phi_{l,k}(v,u)\rho_{t,l}^{\infty}(v)dv\\
&\text{\quad\quad}+\rho^\infty_{t,k-1}(u)\sum_{l=1}^\infty\int_{\mathbb{T}}\phi_{l,k-1}(v,u)\rho^\infty_{t,l}(v)dv
+\rho_{t,k+1}^\infty(u)\sum_{l=0}^\infty\int_{\mathbb{T}}\phi_{k+1,l}(u,v)\rho_{t,l}^\infty(v)dv\\
&\text{\quad\quad if~}k\geq 1,\\
&\frac{d}{dt}\rho_{t,0}^\infty(u)=\rho_{t,1}^{\infty}(u)\sum_{l=0}^\infty\int_{\mathbb{T}}\phi_{1,l}(u,v)\rho_{t,l}^\infty(v)dv
-\rho_{t,0}^\infty(u)\sum_{l=1}^\infty\int_{\mathbb{T}}\phi_{l,0}(v,u)\rho_{t,l}^\infty(v)dv,\\
&\rho_{0,k}^\infty(u)=e^{-\psi(u)}\frac{\psi(u)^k}{k!}
\text{~for~}k\geq 0.
\end{cases}
\end{equation}
To make the solution to Equation \eqref{equ 2.3 mean field ODE infinite types} unique, we further require that the solution to Equation \eqref{equ 2.3 mean field ODE infinite types} satisfies constraints
\begin{equation}\label{equ constraints one}
\inf_{0\leq t\leq T, u\in \mathbb{T}, \atop k\geq 0}\rho^\infty_{t,k}(u)\geq 0,
\text{~}\sup_{0\leq t\leq T, u\in \mathbb{T}}\left|\sum_{l=0}^{+\infty}\rho_{t,l}^\infty(u)-1\right|=0
\text{~and~}\sup_{0\leq t\leq T, u\in \mathbb{T}}\sum_{l=0}^{+\infty}l^r\rho_{t,l}^\infty(u)<+\infty
\end{equation}
for any $T>0, r=1,2$ and
\begin{equation}\label{equ constraints two}
\limsup_{M\rightarrow+\infty}\frac{1}{M}\log \sup_{u\in \mathbb{T}, 0\leq t\leq T}\left(\sum_{l=M}^{+\infty}l\rho_{t,l}^\infty(u)\right)=-\infty
\end{equation}
for any $T>0$.

Now we give the hydrodynamic limit of our model in the case where $\mathcal{K}=+\infty$.

\begin{theorem}\label{theorem 2.2 mean field limit infinite}
If $\mathcal{K}=+\infty$, then, under Assumption (B), the solution $\{\rho_{t,k}^{\infty}\}_{k\geq 0}$ to Equation \eqref{equ 2.3 mean field ODE infinite types} under constraints \eqref{equ constraints one} and \eqref{equ constraints two} exists for $t\in [0, +\infty)$ and is unique. Furthermore, if $\mathcal{K}=+\infty$, then, under Assumptions (A) and (B),
\begin{equation}\label{equ 2.4 convergence in probability infinite}
\lim_{N\rightarrow+\infty}\mu_{t,k}^N(f)=\int_{\mathbb{T}}\rho_{t,k}^\infty(u)f(u)du
\end{equation}
in probability for any $t\geq 0$, $k\geq 0$ and $f\in C(\mathbb{T})$.
\end{theorem}

Theorems \ref{theorem 2.1 mean field limit finite} and \ref{theorem 2.2 mean field limit infinite} show that the hydrodynamic limit of our model is driven by a nonlinear $\left(C(\mathbb{T})\right)^{\mathcal{K}+1}$-valued ordinary differential equation, which is consistent with a heuristic mean field analysis. In detail, according to the generator $\mathcal{L}_N$ of our model given in Equation \eqref{equ 1.1 generator} and Kolmogorov-Chapman equation,
\begin{align*}
\frac{d}{dt}\mathbb{P}\left(\eta_t(i)=k\right)=&-\frac{1}{N}\sum_{j\neq i}\sum_{l=0}^{\mathcal{K}}\phi_{k,l}\left(i/N, j/N\right)\mathbb{P}\left(\eta_t(i)=k,\eta_t(j)=l\right)\\
&-\frac{1}{N}\sum_{j\neq i}\sum_{l=1}^{\mathcal{K}}\phi_{l,k}\left(j/N, i/N\right)\mathbb{P}\left(\eta_t(i)=k,\eta_t(j)=l\right)\\
&+\frac{1}{N}\sum_{j\neq i}\sum_{l=1}^{\mathcal{K}}\phi_{l,k-1}\left(j/N,i/N\right)\mathbb{P}\left(\eta_t(i)=k-1, \eta_t(j)=l\right)\\
&+\frac{1}{N}\sum_{j\neq i}\sum_{l=0}^{\mathcal{K}}\phi_{k+1,l}\left(i/N, j/N\right)\mathbb{P}\left(\eta_t(i)=k+1, \eta_t(j)=l\right).
\end{align*}
When $N$ is large, for $i\neq j$, $\eta_t(i)$ and $\eta_t(j)$ are intuitively considered to be approximately independent and then $\mathbb{P}\left(\eta_t(i)=k,\eta_t(j)=l\right)\approx \mathbb{P}(\eta_t(i)=k)\mathbb{P}(\eta_t(j)=l)$. Taking $\mathbb{P}\left(\eta_t(i)=k\right)$ as the value of a density function at the position $i/N$, then, using above analysis and the fact that $\frac{1}{N}\sum_{j=1}^Ng(j/N)\rightarrow \int_{\mathbb{T}}g(u)du$,  $\mathbb{P}\left(\eta_t(\lfloor Nu\rfloor)=k\right)$ is approximately driven by $\rho_{t,k}^\mathcal{K}(u)$ and
\[
\mathbb{E}\mu_{t,k}^N(f)\approx \frac{1}{N}\sum_{i=1}^N\rho_{t,k}^\mathcal{K}(i/N)f(i/N)\approx \int_{\mathbb{T}}\rho_{t,k}^\mathcal{K}(u)f(u)du.
\]

We further investigate the central limit theorem corresponding to the above hydrodynamic limit, which is also called as the fluctuation, in the case where $\mathcal{K}<+\infty$. To give our result, we first introduce some notations and definitions. For any integer $m\geq 1$, we denote by $\left(C^\infty\left(\mathbb{T}^m\right)\right)^\prime$ the dual of $C^\infty\left(\mathbb{T}^m\right)$ endowed with the weak topology, i.e., $\nu_n\rightarrow \nu$ in $\left(C^\infty\left(\mathbb{T}^m\right)\right)^\prime$ when and only when
\[
\lim_{n\rightarrow+\infty}\nu_n(f)=\nu(f)
\]
for any $f\in C^\infty(\mathbb{T}^m)$. For any integer $m,l\geq 1$ and any linear operator $\mathcal{P}$ from $C^\infty(\mathbb{T}^m)$ to $C^\infty(\mathbb{T}^l)$, we denote by $\mathcal{P}^*$ the conjugate of $\mathcal{P}$ from $\left(C^\infty\left(\mathbb{T}^l\right)\right)^\prime$ to $\left(C^\infty\left(\mathbb{T}^m\right)\right)^\prime$, i.e.,
\[
\mathcal{P}^*\nu(f)=\nu(\mathcal{P}f)
\]
for any $\nu\in \left(C^\infty\left(\mathbb{T}^l\right)\right)^\prime$ and $f\in C^\infty(\mathbb{T}^m)$. We denote by $\{\mathcal{W}_t\}_{t\geq 0}$ the $(C^\infty(\mathbb{T}^2))^\prime$-valued standard Brownian motion, i.e., for any $f\in C^\infty(\mathbb{T}^2)$, $\{\mathcal{W}_t(f)\}_{t\geq 0}$ is a real-valued Brownian motion such that
\[
{\rm Cov}\left(\mathcal{W}_t(f), \mathcal{W}_t(f)\right)=t\int_{\mathbb{T}^2}f^2(u,v)dudv
\]
for all $t\geq 0$. When $\mathcal{K}<+\infty$, we define $\{\mathcal{W}_{t,m,l}: t\geq 0\}_{0\leq m,l\leq \mathcal{K}}$ as a sequence of independent copies of $\{\mathcal{W}_t\}_{t\geq 0}$. For $0\leq k,m,l\leq \mathcal{K}$ and $t\geq 0$, we define $b_{t,m,l}^k$ as the linear operator from $C^\infty(\mathbb{T})$ to $C^\infty\left(\mathbb{T}^2\right)$ such that
\begin{align*}
&b_{t,m,l}^kf(u,v)=\\
&\begin{cases}
\sqrt{\phi_{k,k-1}(u,v)\rho_{t,k}^{\mathcal{K}}(u)\rho_{t,k-1}^{\mathcal{K}}(v)}(f(v)-f(u)) & \text{~if~}m=k\text{~and~}l=k-1,\\
-\sqrt{\phi_{k,k}(u,v)\rho_{t,k}^{\mathcal{K}}(u)\rho_{t,k}^{\mathcal{K}}(v)}(f(v)+f(u)) &  \text{~if~}m=k\text{~and~}l=k, \\
-\sqrt{\phi_{k,l}(u,v)\rho_{t,k}^{\mathcal{K}}(u)\rho_{t,l}^{\mathcal{K}}(v)}(f(u)) & \text{~if~}m=k\text{~and~}l\neq k, k-1,\\
\sqrt{\phi_{k+1,k}(u,v)\rho_{t,k+1}^{\mathcal{K}}(u)\rho_{t,k}^{\mathcal{K}}(v)}(f(u)-f(v)) & \text{~if~}m=k+1\text{~and~}l=k,\\
-\sqrt{\phi_{m,k}(u,v)\rho_{t,m}^{\mathcal{K}}(u)\rho_{t,k}^{\mathcal{K}}(v)}(f(v)) & \text{~if~}m\neq k+1,k\text{~and~}l=k,\\
\sqrt{\phi_{k+1,k-1}(u,v)\rho_{t,k+1}^{\mathcal{K}}(u)\rho_{t,k-1}^{\mathcal{K}}(v)}(f(u)+f(v)) & \text{~if~}m=k+1 \text{~and~}l=k-1,\\
\sqrt{\phi_{m,k-1}(u,v)\rho_{t,m}^{\mathcal{K}}(u)\rho_{t,k-1}^{\mathcal{K}}(v)}(f(v)) & \text{~if~}m\neq k+1,k \text{~and~}l=k-1,\\
\sqrt{\phi_{k+1,l}(u,v)\rho_{t,k+1}^{\mathcal{K}}(u)\rho_{t,l}^{\mathcal{K}}(v)}(f(u)) & \text{~if~}m=k+1 \text{~and~}l\neq k, k-1,\\
0 & \text{~else}
\end{cases}
\end{align*}
for any $f\in C^\infty(\mathbb{T})$ and $u,v\in \mathbb{T}$, where $\{\rho_{t,k}^\mathcal{K}\}_{0\leq k\leq \mathcal{K}}$ is the unique solution to Equation \eqref{equ 2.1 mean field ODE finite types}. Roughly speaking, $(b_{t,m,l}^k)^*d\mathcal{W}_{t,m,l}$ is the fluctuation with respect to $\mu_{t,k}^N$ caused by a particle jumping from a position with $m$ particles to another position with $l$ particles. For $0\leq k,m,l\leq \mathcal{K}$, $t\geq 0$ and $r=1,2$, we define $\mathcal{P}_{t,m,l}^{k,r}$ as the linear operator from $C^\infty(\mathbb{T})$ to $C^\infty(\mathbb{T})$ such that
\begin{align*}
&\mathcal{P}_{t,m,l}^{k,r}f(u)=\\
&\begin{cases}
\int_{\mathbb{T}}\phi_{k,k-1}(u,v)\rho_{t,k-1}^\mathcal{K}(v)(f(v)-f(u))dv & \text{~if~}m=k, l=k-1 \text{~and~}r=1,\\
\int_{\mathbb{T}}\phi_{k,k-1}(v,u)\rho_{t,k}^\mathcal{K}(v)(f(u)-f(v))dv & \text{~if~}m=k, l=k-1 \text{~and~}r=2,\\
-\int_{\mathbb{T}}\phi_{k,k}(u,v)\rho_{t,k}^\mathcal{K}(v)(f(u)+f(v))dv & \text{~if~}m=k, l=k \text{~and~}r=1,\\
-\int_{\mathbb{T}}\phi_{k,k}(v,u)\rho_{t,k}^\mathcal{K}(v)(f(u)+f(v))dv & \text{~if~}m=k, l=k \text{~and~}r=2,\\
-\int_{\mathbb{T}}\phi_{k,l}(u,v)\rho_{t,l}^\mathcal{K}(v)f(u)dv & \text{~if~}m=k, l\neq k, k-1 \text{~and~}r=1,\\
-\int_{\mathbb{T}}\phi_{k,l}(v,u)\rho_{t,k}^\mathcal{K}(v)f(v)dv & \text{~if~}m=k, l\neq k, k-1 \text{~and~}r=2,\\
\int_{\mathbb{T}}\phi_{k+1,k}(u,v)\rho_{t,k}^\mathcal{K}(v)(f(u)-f(v))dv & \text{~if~}m=k+1, l=k \text{~and~}r=1,\\
\int_{\mathbb{T}}\phi_{k+1,k}(v,u)\rho_{t,k+1}^\mathcal{K}(v)(f(v)-f(u))dv & \text{~if~}m=k+1, l=k \text{~and~}r=2,\\
-\int_{\mathbb{T}}\phi_{m,k}(u,v)\rho_{t,k}^\mathcal{K}(v)f(v)dv & \text{~if~}m\neq k+1,k, l=k \text{~and~}r=1,\\
-\int_{\mathbb{T}}\phi_{m,k}(v,u)\rho_{t,m}^\mathcal{K}(v)f(u)dv & \text{~if~}m\neq k+1,k, l=k \text{~and~}r=2,\\
\int_{\mathbb{T}}\phi_{k+1,k-1}(u,v)\rho_{t,k-1}^\mathcal{K}(v)(f(u)+f(v))dv & \text{~if~}m=k+1, l=k-1 \text{~and~}r=1,\\
\int_{\mathbb{T}}\phi_{k+1,k-1}(v,u)\rho_{t,k+1}^\mathcal{K}(v)(f(u)+f(v))dv & \text{~if~}m=k+1, l=k-1 \text{~and~}r=2,\\
\int_{\mathbb{T}}\phi_{m,k-1}(u,v)\rho_{t,k-1}^\mathcal{K}(v)f(v)dv & \text{~if~}m\neq k+1,k, l=k-1 \text{~and~}r=1,\\
\int_{\mathbb{T}}\phi_{m,k-1}(v,u)\rho_{t,m}^\mathcal{K}(v)f(u)dv & \text{~if~}m\neq k+1,k, l=k-1 \text{~and~}r=2,\\
\int_{\mathbb{T}}\phi_{k+1,l}(u,v)\rho_{t,l}^\mathcal{K}(v)f(u)dv & \text{~if~}m=k+1, l\neq k,k-1 \text{~and~}r=1,\\
\int_{\mathbb{T}}\phi_{k+1,l}(v,u)\rho_{t,k+1}^\mathcal{K}(v)f(v)dv & \text{~if~}m=k+1, l\neq k,k-1 \text{~and~}r=2,\\
0 & \text{~else}
\end{cases}
\end{align*}
for any $f\in C^\infty(\mathbb{T})$ and $u\in \mathbb{T}$. In this paper, when $\mathcal{K}<+\infty$, we are concerned with the following $\left(\left(C^\infty(\mathbb{T})\right)^\prime\right)^{\mathcal{K}+1}$-valued O-U process $\{V_{t,k}: t\geq 0\}_{0\leq k\leq \mathcal{K}}$: for all $0\leq k\leq \mathcal{K}$,
\begin{equation}\label{equ O-U}
dV_{t,k}=\sum_{m=0}^\mathcal{K}\sum_{l=0}^\mathcal{K}\left(\left(\mathcal{P}_{t,m,l}^{k,1}\right)^*dV_{t,m}
+\left(\mathcal{P}_{t,m,l}^{k,2}\right)^*dV_{t,l}
+(b_{t,m,l}^k)^*d\mathcal{W}_{t,m,l}\right).
\end{equation}
The solution $\{V_{t,k}: t\geq 0\}_{0\leq k\leq \mathcal{K}}$ to Equation \eqref{equ O-U} is rigorously defined as the solution to the corresponding martingale problem. In detail, using an analysis similar with that in the proof of Theorem (1.4) of \cite{Holley1978}, when the initial distribution of $\{V_{0,k}\}_{0\leq k\leq \mathcal{K}}$ is given, there exists a unique random element $\{V_{t,k}: t\geq 0\}_{0\leq k\leq \mathcal{K}}$ in $C\left([0, +\infty), \left(\left(C^\infty(\mathbb{T})\right)^\prime\right)^{\mathcal{K}+1}\right)$ such that
\begin{align*}
\Bigg\{
G&\left(\{V_{t,k}(f_k)\}_{k=0}^{\mathcal{K}}\right)-G\left(\{V_{0,k}(f_k)\}_{k=0}^{\mathcal{K}}\right)\\
&-\sum_{k=0}^\mathcal{K}\sum_{m=0}^\mathcal{K}\sum_{l=0}^\mathcal{K}\int_0^t\partial_kG\left(\{V_{s,k}(f_k)\}_{k=0}^{\mathcal{K}}\right)
\left(V_{s,m}(\mathcal{P}_{t,m,l}^{k,1}f_k)+V_{s,l}(\mathcal{P}_{t,m,2}^{k,1}f_k)\right)ds\\
&-\frac{1}{2}\sum_{k_1=0}^\mathcal{K}\sum_{k_2=0}^\mathcal{K}\sum_{m=0}^\mathcal{K}\sum_{l=0}^\mathcal{K}\int_0^t\int_{\mathbb{T}^2}\partial^2_{k_1k_2}
G\left(\{V_{0,k}(f_k)\}_{k=0}^{\mathcal{K}}\right)(b_{s,m,l}^{k_1}f_{k_1}(u,v))(b_{s,m,l}^{k_2}f_{k_2}(u,v))dsdudv\Bigg\}_{t\geq 0}
\end{align*}
is a martingale for any $f_0, f_1,\ldots, f_\mathcal{K}\in C^\infty(\mathbb{T})$ and $G\in C_c^\infty(\mathbb{R}^{\mathcal{K}+1})$, where $\partial_k$ is the partial derivative with respect to the $k$-th coordinate. This martingale-problem solution $\{V_{t,k}: t\geq 0\}_{0\leq k\leq \mathcal{K}}$ is defined as the solution to Equation \eqref{equ O-U} under a given initial distribution.

Now we define our fluctuation density fields. For any $N\geq 1$, $t\geq 0$ and $0\leq k\leq \mathcal{K}$, the fluctuation density field $V_{t,k}^N$ with
respect to $\mu_{t,k}^N$ is defined as
\[
V_{t,k}^N(du)=\frac{1}{\sqrt{N}}\sum_{i=1}^N\left(1_{\{\eta_t^N(i)=k\}}-\mathbb{P}(\eta_t^N(i)=k)\right)\delta_{i/N}(du).
\]
Let $\left(C^\infty(\mathbb{T})\right)^\prime$ be endowed with the weak topology, under which $\lim_{n\rightarrow+\infty}\nu_n=\nu$ in $\left(C^\infty(\mathbb{T})\right)^\prime$ when and only when $\lim_{n\rightarrow+\infty}\nu_n(f)=\nu(f)$ for all $f\in C^\infty(\mathbb{T})$. Then, for given $T>0$, $\{V_{t,k}^N: 0\leq t\leq T\}_{k=0}^\mathcal{K}$ is a random element in the set of c\`{a}dl\`{a}g functions from $[0, T]$ to $\left(\left(C^\infty(\mathbb{T})\right)^\prime\right)^{\mathcal{K}+1}$, which we denote by $D\left([0, T], \left(\left(C^\infty(\mathbb{T})\right)^\prime\right)^{\mathcal{K}+1}\right)$. We denote by $D([0, T], \mathbb{R})$ the set of c\`{a}dl\`{a}g functions from $[0, T]$ to $\mathbb{R}$ and let $D([0, T], \mathbb{R})$ be endowed with the Skorohod typology. Let $D\left([0, T], \left(\left(C^\infty(\mathbb{T})\right)^\prime\right)^{\mathcal{K}+1}\right)$ be endowed with the topology under which $\lim_{n\rightarrow+\infty}\theta_n=\theta$ in $D\left([0, T], \left(\left(C^\infty(\mathbb{T})\right)^\prime\right)^{\mathcal{K}+1}\right)$ when and only when
$\lim_{n\rightarrow+\infty}\{\theta_{n,t,k}(f)\}_{0\leq t\leq T}=\{\theta_{t,k}(f)\}_{0\leq t\leq T}$ in $D([0, T], \mathbb{R})$ for any $f\in C^\infty(\mathbb{T})$ and $0\leq k\leq \mathcal{K}$.

Now we give the third  main result of this paper, which is about the fluctuation limit of our model in the case where $\mathcal{K}<+\infty$.

\begin{theorem}\label{theorem 2.3 fluctuation limit}
If $\mathcal{K}<+\infty$, then, under Assumption (A), $\{V_{t,k}^N: 0\leq t\leq T\}_{k=0}^{\mathcal{K}}$ converges weakly, with respect to the topology of $D\left([0, T], \left(\left(C^\infty(\mathbb{T})\right)^\prime\right)^{\mathcal{K}+1}\right)$, to $\{V_{t,k}: 0\leq t\leq T \}_{k=0}^{\mathcal{K}}$, where $\{V_{t,k}: 0\leq t\leq T \}_{k=0}^{\mathcal{K}}$ is the unique solution to Equation \eqref{equ O-U} with the Gaussian initial distribution under which
\[
\mathbb{E}V_{0,k}(f)=0
\]
and
\begin{align*}
&{\rm Cov}\left(V_{0,k}(f), V_{0, m}(g)\right)=\\
&\begin{cases}
&-\int_{\mathbb{T}}{\mathcal{K} \choose k}{\mathcal{K} \choose m}\left(\psi(u)/\mathcal{K}\right)^{k+m}\left(1-\frac{\psi(u)}{\mathcal{K}}\right)^{2\mathcal{K}-k-m}f(u)g(u)du  \text{\quad if~}m\neq k,\\
&\int_{\mathbb{T}}{\mathcal{K} \choose k}\left(\psi(u)/\mathcal{K}\right)^k\left(1-\frac{\psi(u)}{\mathcal{K}}\right)^{\mathcal{K}-k}
\left(1-{\mathcal{K} \choose k}\left(\psi(u)/\mathcal{K}\right)^k\left(1-\frac{\psi(u)}{\mathcal{K}}\right)^{\mathcal{K}-k}\right)f(u)g(u)du  \\
&\text{\quad\quad\quad}\text{\quad if~}m=k
\end{cases}
\end{align*}
for any $0\leq k, m\leq \mathcal{K}$ and $f,g\in C^\infty(\mathbb{T})$.
\end{theorem}

A heuristic explanation of Theorem \ref{theorem 2.3 fluctuation limit} is as follows. Using Dynkin's martingale formula,
\[
\Big\{G\left(\{V_{t,k}^N(f_k)\}_{k=0}^{\mathcal{K}}\right)-G\left(\{V_{0,k}^N(f_k)\}_{k=0}^{\mathcal{K}}\right)
-\int_0^t(\partial_s+\mathcal{L}_N)G\left(\{V_{s,k}^N(f_k)\}_{k=0}^{\mathcal{K}}\right)ds\Big\}_{t\geq 0}
\]
is a martingale. Using Taylor's expansion formula on $\mathcal{L}_NG\left(\{V_{s,k}^N(f_k)\}_{k=0}^{\mathcal{K}}\right)$ up to the second order, we obtain an approximated expression of $(\partial_s+\mathcal{L}_N)G\left(\{V_{s,k}^N(f_k)\}_{k=0}^{\mathcal{K}}\right)$. In this approximation, with small errors, we can further replace a term with the form
\[
\frac{1}{N}\sum_{i=1}^N1_{\{\eta_t^N(i)=k\}}f(i/N) \text{~or~} \frac{1}{N}\sum_{i=1}^N\mathbb{P}(\eta_t^N(i)=k)f(i/N)
\]
by
\[
\mu_{t,k}(f)=\int_{\mathbb{T}}\rho^{\mathcal{K}}_{t,k}(u)f(u)du
\]
according to our mean field limit result. Furthermore, using the approximated independence of numbers of particles on different positions, we can replace
$\mathbb{P}(\eta_t^N(i)=k, \eta_t^N(j)=m)$ by $\mathbb{P}(\eta_t^N(i)=k)\mathbb{P}(\eta_t^N(j)=m)$ for $i\neq j$ and then replace a term with the form
\[
\frac{1}{N^{3/2}}\sum_{i=1}^N\sum_{j=1}^N\left(1_{\{\eta_t^N(i)=k, \eta_t^N(j)=m\}}-\mathbb{P}(\eta_t^N(i)=k, \eta_t^N(j)=m)\right)g(i/N, j/N)
\]
by
\begin{align*}
&\frac{1}{N^{1/2}}\sum_{i=1}^N\left(1_{\{\eta_t^N(i)=k\}}-\mathbb{P}(\eta_t^N(i)=k)\right)\left(\int_{\mathbb{T}}g(i/N,v)\rho_{t,m}^\mathcal{K}(v)dv\right)\\
&+\frac{1}{N^{1/2}}\sum_{j=1}^N\left(1_{\{\eta_t^N(j)=m\}}-\mathbb{P}(\eta_t^N(j)=m)\right)\left(\int_{\mathbb{T}}g(u,j/N)\rho_{t,k}^\mathcal{K}(u)du\right)\\
&=V_{t,k}^N\left(\int_{\mathbb{T}}g(\cdot, v)\rho_{t,m}^{\mathcal{K}}(v)dv\right)+V_{t,m}^N\left(\int_\mathbb{T}g(v,\cdot)\rho_{t,k}^\mathcal{K}(v)dv\right)
\end{align*}
with small errors. Using above approximations, $\{V_{t,k}^N: t\geq 0\}_{k=0}^\mathcal{K}$ is an approximated solution to the martingale problem with respect to Equation \eqref{equ O-U}, which intuitively implies Theorem \ref{theorem 2.3 fluctuation limit}.

It is natural to ask whether an analogue fluctuation theorem of Theorem \ref{theorem 2.3 fluctuation limit} holds in the case where $\mathcal{K}=+\infty$. It is natural to guess that the conclusion is positive and the corresponding fluctuation equation should be Equation \eqref{equ O-U} with the finite $\mathcal{K}$ replaced by $+\infty$. To prove such a fluctuation theorem, we need to show that the covariance of numbers of particles on two different positions is with order $O(N^{-1})$ according to our above strategy. However, currently we can only give this covariance estimation in the case where $\mathcal{K}<\infty$. In the case where $\mathcal{K}=+\infty$, we can only show that the aforesaid covariance is with order $O(N^{-0.9})$, which is enough for the proof of the hydrodynamic limit but not enough for the fluctuation. We will work on the fluctuation in the case where $\mathcal{K}=+\infty$ as a further investigation.

The rest of this paper is arranged as follows. In Section \ref{section three}, using a graphical method, we prove the approximated independence of the numbers of particles on different positions, which plays key roles in proofs of our main results as we have shown in heuristic explanations of Theorems \ref{theorem 2.1 mean field limit finite}-\ref{theorem 2.3 fluctuation limit}. In Sections \ref{section four proof of hydro finite} and \ref{section five proof of hydro infinite}, we prove Theorems \ref{theorem 2.1 mean field limit finite} and \ref{theorem 2.2 mean field limit infinite}. Using the approximated independence given in Section \ref{section three}, these two theorems follow from the Kolmogorov-Chapman equation of our model. In Section \ref{section six proof of fluctuation}, we prove Theorem \ref{theorem 2.3 fluctuation limit}. Our proof utilizes a Dynkin's martingale strategy, the outline of which has been given in the heuristic explanation of Theorem \ref{theorem 2.3 fluctuation limit}. In Section \ref{section seven}, we apply our main theorems on examples such as exclusion processes and Ehrenfest models.

\section{Approximated independence}\label{section three}

In this section we will prove following two lemmas.

\begin{lemma}\label{lemma 3.1 approximated independence finite K}
Let $T>0$ be a given moment. If $\mathcal{K}<+\infty$, then, under Assumption (A), there exists $C_2<+\infty$ independent of $N$ such that
\begin{equation}\label{equ finite K two positions approximate independence}
\left|\mathbb{P}\left(\eta_t^N(x)=k_1, \eta_t^N(y)=k_2\right)-\mathbb{P}(\eta_t^N(x)=k_1)\mathbb{P}(\eta_t^N(y)=k_2)\right|\leq \frac{C_2}{N}
\end{equation}
for all $0\leq t\leq T$, $N\geq 1$, $k_1, k_2\in \{0,1,\ldots,\mathcal{K}\}$ and any two different $x,y$ in $\{1,2,\ldots,N\}$.
\end{lemma}

\begin{lemma}\label{lemma 3.2 approximated independence infinite K}
Let $T>0$ be a given moment. If $\mathcal{K}=+\infty$, then, under Assumptions (A) and (B), there exists an integer $N_0\geq 1$ such that
\begin{equation}\label{equ infinite K two posistitons approximate independence}
\left|\mathbb{P}\left(\eta_t^N(x)=k_1, \eta_t^N(y)=k_2\right)-\mathbb{P}(\eta_t^N(x)=k_1)\mathbb{P}(\eta_t^N(y)=k_2)\right|\leq \frac{1}{N^{0.9}}
\end{equation}
for all $N\geq N_0$, $0\leq t\leq T$, $k_1, k_2\geq 0$ and any two different $x,y\in \{1,2,\ldots,N\}$.
\end{lemma}

We prove Lemmas \ref{lemma 3.1 approximated independence finite K} and \ref{lemma 3.2 approximated independence infinite K} in Subsections \ref{subsection 3.1} and \ref{subsection 3.2} respectively. Proofs of these two lemmas follow similar strategies. For any two different $x,y\in \{1,2,\ldots,N\}$, we use a graphical method to construct independent $\widetilde{\eta}_t^N(x), \widetilde{\eta}_t^N(y)$ such that $\mathbb{P}(\widetilde{\eta}_{t}^N(z)\neq \eta_t^N(z))=O(N^{-1})$ (resp. $O(N^{-0.9})$) for $z=x,y$ in the case where $\mathcal{K}<+\infty$ (resp. $\mathcal{K}=+\infty$). The difference between covariance-estimations in above two cases, according to our current approach, is mainly due to the fact that the upper bound $\sup_{k,l,u,v}\phi_{k,l}(u,v)$ of jumping rates is finite in the case where $\mathcal{K}<+\infty$ but might be infinite in the case where $\mathcal{K}=+\infty$.
As we have explained in Section \ref{section two}, if we can improve $N^{-0.9}$ in the right-hand side of \eqref{equ infinite K two posistitons approximate independence} to $N^{-1}$, then we can prove an analogue of Theorem \ref{theorem 2.3 fluctuation limit} in the case where $\mathcal{K}=+\infty$. We will work on this improvement as a further investigation.

\subsection{The proof of Lemma \ref{lemma 3.1 approximated independence finite K}}\label{subsection 3.1}

In this subsection, we prove Lemma \ref{lemma 3.1 approximated independence finite K}. Throughout this subsection, we assume that $\mathcal{K}<+\infty$. We first utilize a graphical method to give an equivalent definition of our process $\{\eta_t^N\}_{t\geq 0}$. Let
\[
K_1=\sup_{0\leq k,l\leq \mathcal{K},u,v\in \mathbb{T}}\phi_{k,l}(u,v).
\]
Note that $K_1<+\infty$ since $\mathcal{K}<+\infty$. For any two different $x,y\in \{1,2,\ldots, N\}$, let $\{\Xi_t^{(x,y)}\}_{t\geq 0}$ be a Poisson process with rate $\frac{K_1}{N}$ and $\{U^{(x,y)}_k\}_{k\geq 1}$ be independent and identically distributed with the uniform distribution on $[0, 1]$. Note that we care the order of $x$ and $y$, i.e., $U^{(x,y)}_k\neq U^{(y,x)}_k$ and $\Xi_t^{(x,y)}\neq \Xi_t^{(y,x)}$. We assume that all above Poisson processes and uniformly distributed random variables are independent. Since the number of above Poisson processes is finite, we can arrange all event moments of these Poisson processes as
\[
t_1<t_2<t_3<\ldots.
\]
For each $k\geq 1$, $t_k$ is the $b_k$-th event moment of $\{\Xi^{(x_k, y_k)}_t\}_{t\geq 0}$ for some $x_k ,y_k\in \{1,2,\ldots,N\}$ and some integer $b_k\geq 1$. Now we define $\{\eta_t^N\}_{t\geq 0}$ by induction. We supplementarily define $t_0=0$. For $t=0$, let $\eta_0^N$ be distributed according to Assumption (A). For $0\leq t<t_1$, let $\eta_t^N=\eta_0^N$. Assuming that $\eta_t^N$ has been defined for $0\leq t<t_k$ for some $k\geq 1$. At the moment $t_k$, if
\[
U^{(x_k, y_k)}_{b_k}\leq \frac{\phi_{\eta^N_{t_{k-1}}(x_k), \eta_{t_{k-1}}^N(y_k)}\left(x_k/N, y_k/N\right)}{K_1},
\]
then we let $\eta_{t_k}^N=\left(\eta_{t_{k-1}}^N\right)^{(x_k,y_k)}$, otherwise we let $\eta_{t_k}^N=\eta_{t_{k-1}}^N$. For $t_k<t<t_{k+1}$, we let $\eta_t^N=\eta_{t_k}^N$ and consequently $\eta_t^N$ is defined for $0\leq t<t_{k+1}$. By induction, $\eta_t^N$ is defined for all $t\geq 0$.

The above construction is called a graphical method in the sense that for each $k\geq 1$ we put an arrow on the product graph $\{1,2,\ldots,N\}\times [0,+\infty)$ from $(x_k,t_k)$ to $(y_k, t_k)$ and a particle jumps through this arrow from $x_k$ to $y_k$ when and only when an independent copy of a random variable $U$ uniformly distributed on $[0, 1]$ is at most $\phi_{m,l}(x_k/N, y_k/N)/K_1$, where $m,l$ are numbers of particles on $x_k$ and $y_k$ respectively at the moment just before $t_k$. When there are $m$ particles on position $x$ and $l$ particles on position $y$, since an arrow from $x$ to $y$ occurs at rate $K_1/N$ and a particle jumps through this arrow with probability
\[
\mathbb{P}(U\leq \phi_{m,l}(x/N, y/N)/K_1)=\phi_{m,l}(x/N, y/N)/K_1,
\]
we have that a particle jumps from $x$ to $y$ at rate
\[
\frac{K_1}{N}\times \phi_{m,l}(x/N, y/N)/K_1=\frac{1}{N}\phi_{m,l}(x/N, y/N).
\]
Consequently, when the initial distribution is given, the process constructed according to the above graphical method is identically distributed with the one defined via the generator given by \eqref{equ 1.1 generator}. For a comprehensive reading of the theory of graphical methods of interacting particle systems, see \cite{Har1978} or Section 3.6 of \cite{Lig1985}.

Now we introduce the definition of an influence path. For any two different $x,y\in \{1,2,\ldots, N\}$ and $s>0$, we write $(x,s)\sim (y,s)$ when and only when $s$ is an event moment of $\{\Xi_t^{(x,y)}\}_{t\geq 0}$ or $\{\Xi_t^{(y,x)}\}_{t\geq 0}$. For any two different $x,y\in \{1,2,\ldots, N\}$ and $t>0$, if there exist an integer $m\geq 1$, $x_1, x_2,\ldots, x_{m-1}\in \{1,2,\ldots,N\}$ and $0<t_1<t_2<\ldots<t_m\leq t$ such that $(x_{k-1}, t_k)\sim (x_k, t_k)$ for all $1\leq k\leq m$, where $x_0=y$ and $x_m=x$, then we say that $y,x_1,\ldots,x_{m-1},x$ is a $t$-influence path. For any $t>0$ and $x\in \{1,2,\ldots, N\}$, we define
\[
\Gamma_{t,x}=\{x\}\cup\{y:\text{~there is a $t$-influence path from $y$ to $x$}\}.
\]
According to the graphical construction of our process, the number of particles on a position $x$ may change only at moments when there is an arrow starting from $x$ or ending at $x$. Consequently, the value of $\eta_t^N(x)$ is determined by $\Gamma_{t,x}$, $\{\Xi_s^{(z,w)}:~0\leq s\leq t\}_{z,w\in \Gamma_{t,x}}$,
$\{U^{(z,w)}_k:~k\geq 1\}_{z,w\in \Gamma_{t,x}}$ and $\{\eta_0^N(z)\}_{z\in \Gamma_{t,x}}$. To give this property a rigorous mathematical description, we denote by $\Lambda_t$ the set of nondecreasing c\`{a}dl\`{a}g functions from $[0, T]$ to $\{0,1,2,\ldots\}$ starting at $0$ and by $\beta_B$ the set
\[
\{(z,w):~z,w\in B, z\neq w\}
\]
for any $B\subseteq \{0,1,2,\ldots,N\}$. For any $x\in \{1,2,\ldots,N\}$ and $B\ni x$, conditioned on $\Gamma_{t,x}=B$, there exists a function $Q_{x, B}$ from
\[
(\Lambda_t)^{\beta_B}\times ([0,1]^{\{1,2,3,\ldots\}})^{\beta_B}\times (\mathbb{Z}_{\mathcal{K}})^B
\]
to $\mathbb{Z}_{\mathcal{K}}$ such that
\begin{equation}\label{equ Q}
\eta_t^N(x)=Q_{x, B}\left(\{\Xi_s^{(z,w)}:~0\leq s\leq t\}_{(z,w)\in \beta_B}, \{U^{(z,w)}_k:~k\geq 1\}_{(z,w)\in \beta_B}, \{\eta_0^N(z)\}_{z\in B}\right).
\end{equation}
According to a breadth-first-search (BFS), $\Gamma_{t,x}=\bigcup_{m=0}^N\Gamma_{t,x}^m$, where $\{\Gamma_{t,x}^m\}_{0\leq m\leq N}$ is defined as follows.

(1) $\Gamma_{t,x}^0=\{x\}$.

(2) $\Gamma_{t,x}^1$ is the set of $y$ such that there exists $0<t_1\leq t$ making $(y, t_1)\sim (x, t_1)$.

(3) Assuming that $\Gamma_{t,x}^r$ has been defined for $1\leq r\leq m$, then $\Gamma_{t,x}^{m+1}$ is the set of $y\not\in \bigcup_{r=0}^m\Gamma_{t,x}^r$ such that there exist $0<t_1<t_2<\ldots<t_{m+1}\leq t$ and $x_1, x_2,\ldots, x_m\in \bigcup_{r=0}^m\Gamma_{t,x}^r$ making $(x_{k-1}, t_k)\sim (x_k, t_k)$ for all $1\leq k\leq m+1$, where $x_0=y$ and $x_{m+1}=x$.

In brief, $\Gamma_{t,x}^m$ is the set of $y$ such that there is a $t$-influence path from $y$ to $x$ with length $m$ but there is no $t$-influence path from $y$ to $x$ with length at most $m-1$.

To give the proof of Lemma \ref{lemma 3.1 approximated independence finite K}, we need further introduce some notations and definitions for later use. For any two different $z,w\in \{0,1,\ldots, N\}$, we define $\{\hat{\Xi}_t^{(z,w)}\}_{t\geq 0}$, $\{\hat{U}^{(z,w)}_k\}_{k\geq 1}$ and $\hat{\eta}_0^N(z)$ as independent copies of $\{\Xi_t^{(z,w)}\}_{t\geq 0}$, $\{U^{(z,w)}_k\}_{k\geq 1}$ and $\eta_0^N(z)$ respectively. We further assume that all these Poisson processes, uniformly distributed random variables and $\{\hat{\eta}_0^N(z)\}_{1\leq z\leq N}$ are independent. For given $x\in \{1,2,\ldots,N\}$ and any two different $z,w\in \{1,2,\ldots, N\}$, we further define
\[
\widetilde{\Xi}_t^{(z,w)}
=
\begin{cases}
\Xi_t^{(z,w)} &\text{~if~}z,w\not\in \Gamma_{t,x},\\
\hat{\Xi}_t^{(z,w)} & \text{~if~}z\in \Gamma_{t,x} \text{~or~}w\in \Gamma_{t,x}
\end{cases}
\]
for all $t\geq 0$,
\[
\widetilde{U}_k^{(z,w)}
=
\begin{cases}
U_k^{(z,w)} &\text{~if~}z,w\not\in \Gamma_{t,x},\\
\hat{U}_k^{(z,w)} & \text{~if~}z\in \Gamma_{t,x} \text{~or~}w\in \Gamma_{t,x}
\end{cases}
\]
for all $k\geq 1$ and
\[
\widetilde{\eta}_0^N(z)
=
\begin{cases}
\eta_0^N(z) &\text{~if~}z\not\in \Gamma_{t,x},\\
\hat{\eta}_0^N(z) & \text{~if~}z\in \Gamma_{t,x}.
\end{cases}
\]
For any $y\neq x$, we define $\widetilde{\Gamma}_{t,y}$ as the analogue of $\Gamma_{t,y}$ generated by $\{\widetilde{\Xi}_t^{(z,w)}:~t\geq 0\}_{z\neq w}$ and $\{\widetilde{U}_k^{(z,w)}:~k\geq 1\}_{z\neq w}$. That is to say, $\widetilde{\Gamma}_{t,y}=\bigcup_{m=0}^{N}\widetilde{\Gamma}_{t,y}^m$, where $\{\widetilde{\Gamma}_{t,y}^m\}_{0\leq m\leq N}$ is defined as follows.

(1) $\widetilde{\Gamma}_{t,y}^0=\{y\}$.

(2) $\widetilde{\Gamma}_{t,y}^1$ is the set of $z$ such that there exists $0<t_1\leq t$ making $t_1$ an event moment of $\{\widetilde{\Xi}_s^{(y,z)}\}_{s\geq 0}$ or $\{\widetilde{\Xi}_s^{(z,y)}\}_{s\geq 0}$.

(3) Assuming that $\widetilde{\Gamma}_{t,y}^r$ has been defined for $1\leq r\leq m$, then $\widetilde{\Gamma}_{t,y}^{m+1}$ is the set of $z\not\in \bigcup_{r=0}^m\widetilde{\Gamma}_{t,y}^r$ such that there exist $0<t_1<t_2<\ldots<t_{m+1}\leq t$ and $z_1, z_2,\ldots, z_m\in \bigcup_{r=0}^m\widetilde{\Gamma}_{t,y}^r$ making $t_k$ an event moment of $\{\widetilde{\Xi}_s^{(z_{k-1}, z_k)}\}_{s\geq 0}$ or $\{\widetilde{\Xi}_s^{(z_k, z_{k-1})}\}_{s\geq 0}$ for all $1\leq k\leq m+1$, where $z_0=z$ and $z_{m+1}=y$.

Now we define
\[
\widetilde{\eta}_t^N(y)=
Q_{y, \widetilde{\Gamma}_{t,y}}\left(\{\widetilde{\Xi}_s^{(z,w)}:~0\leq s\leq t\}_{(z,w)\in \beta_{\widetilde{\Gamma}_{t,y}}}, \{\widetilde{U}^{(z,w)}_k:~k\geq 1\}_{(z,w)\in \beta_{\widetilde{\Gamma}_{t,y}}}, \{\widetilde{\eta}_0^N(z)\}_{z\in \widetilde{\Gamma}_{t,y}}\right),
\]
where $\{Q_{x, B}\}_{1\leq x\leq N, B\ni x}$ is defined as in \eqref{equ Q}. The following three lemmas are crucial for us to prove Lemma \ref{lemma 3.1 approximated independence finite K}.

\begin{lemma}\label{lemma 3.1.1}
For given $x\in \{1,2,\ldots,N\}$ and any subset $B$ of $\{1,2,\ldots, N\}$ such that $x\in B$, the event $\{\Gamma_{t,x}=B\}$ is measurable with the $\sigma$-algebra generated by $\{\Xi_t^{(w,z)}:~t\geq 0, w\in B\text{~or~}z\in B\}$, i.e.,  $1_{\{\Gamma_{t,x}=B\}}$ is independent of $\{\Xi_t^{(w,z)}:~t\geq 0, w,z\not\in B\}$.
\end{lemma}

\begin{lemma}\label{lemma 3.1.2}
For $y\neq x$, $\widetilde{\eta}_t^N(y)$ is independent of $\eta_t^N(x)$.
\end{lemma}

\begin{lemma}\label{lemma 3.1.3}
For given $T>0$, there exists $C_3<+\infty$ independent of $N$ such that
\[
\mathbb{P}\left(\Gamma_{t,x}\bigcap \Gamma_{t,y}\neq \emptyset\right)\leq \frac{C_3}{N}
\]
for all $0\leq t\leq T$, $N\geq 1$ and any two different $x,y\in \{1,2,\ldots,N\}$.
\end{lemma}

\proof[Proof of Lemma \ref{lemma 3.1.1}]

Assuming that a sample of $\{\Xi_t^{(w,z)}:~t\geq 0, w\neq z\}$ makes $\Gamma_{t,x}=B$, then we only need to show that any change of the value of $\{\Xi_t^{(w,z)}\}_{t\geq 0}$ will not make $\Gamma_{t,x}\neq B$ when $w, z\not \in B$. For given $w,z\not\in B$, if a change of $\{\Xi_t^{(w,z)}\}_{t\geq 0}$ makes $\Gamma_{t,x}\neq B$, then there exists an $t$-influence path from $u$ to $x$ utilizing the changed values of Poisson processes for some $u\not\in B$ and this influence path must visit the edge connecting $w,z$. Then, from the last time when this path visits $\{w,z\}$, there is an $t$-influence path from $w$ to $x$ avoiding $z$ or from $z$ to $x$ avoiding $w$, which does not utilize the changed value of $\{\Xi_t^{(w,z)}\}_{t\geq 0}$ and is contradictory to $w,z\not\in \Gamma_{t,x}$ before $\{\Xi_t^{(w,z)}\}_{t\geq 0}$ being changed. Consequently, the proof is complete.

\qed

\proof[Proof of Lemma \ref{lemma 3.1.2}]

We only need to show that, for any given $0\leq k_1, k_2\leq \mathcal{K}$ and any $B\ni x$, there exists $A=A(k_1)$ independent of $k_2$ and $B$ such that
\begin{equation}\label{equ 3.1}
\mathbb{P}\left(\widetilde{\eta}_t^N(y)=k_1, \eta_t^N(x)=k_2\Big|\Gamma_{t,x}=B\right)=A\mathbb{P}\left(\eta_t^N(x)=k_2\Big|\Gamma_{t,x}=B\right).
\end{equation}
According to the definition of $\widetilde{\eta}_t^N(y)$, there exists a deterministic function $\widetilde{Q}_y$ from
\[
\{\Lambda_t\}^{\beta_{\{1,2,\ldots,N\}}}\times ([0,1]^{\{1,2,3,\ldots\}})^{\beta_{\{1,2,\ldots,N\}}}\times (\mathbb{Z}_{\mathcal{K}})^N
\]
such that
\[
\widetilde{\eta}_t^N(y)
=\widetilde{Q}_y\left(\{\widetilde{\Xi}_t^{(z,w)}:~t\geq 0\}_{z\neq w}, \{\widetilde{U}_k^{(z,w)}:~k\geq 1\}_{z\neq w}, \{\widetilde{\eta}_0^N(z)\}_{1\leq z\leq N}\right).
\]
Then, conditioned on $\Gamma_{t,x}=B$,
\begin{align*}
&\widetilde{\eta}_t^N(y)=\widetilde{Q}_y\Bigg(\{\hat{\Xi}_t^{(z,w)}:~t\geq 0\}_{z\in B\text{~or~}w\in B, z\neq w}, \{\Xi_t^{(z,w)}:~t\geq 0\}_{z,w\not\in B, z\neq w},  \\
&\text{\quad}\{\hat{U}_k^{(z,w)}:~k\geq 1\}_{z\in B\text{~or~}w\in B, z\neq w}, \{U_k^{(z,w)}:~k\geq 1\}_{z,w\not\in B, z\neq w}, \{\hat{\eta}_0^N(z)\}_{z\in B}, \{\eta_0^N(z)\}_{z\not\in B}\Bigg)
\end{align*}
and
\[
\eta_t^N(x)=Q_{x, B}\left(\{\Xi_s^{(z,w)}:~0\leq s\leq t\}_{(z,w)\in \beta_B}, \{U^{(z,w)}_k:~k\geq 1\}_{(z,w)\in \beta_B}, \{\eta_0^N(z)\}_{z\in B} \right).
\]
According to Lemma \ref{lemma 3.1.1} and our definitions of $\hat{\Xi}_t^{(z,w)}, \hat{U}_k^{(z,w)}$ and $\hat{\eta}_0^N(z)$, we have that
\begin{align*}
&\widetilde{Q}_y\Bigg(\{\hat{\Xi}_t^{(z,w)}:~t\geq 0\}_{z\in B\text{~or~}w\in B, z\neq w}, \{\Xi_t^{(z,w)}:~t\geq 0\}_{z,w\not\in B, z\neq w},  \\
&\text{\quad}\{\hat{U}_k^{(z,w)}:~k\geq 1\}_{z\in B\text{~or~}w\in B, z\neq w}, \{U_k^{(z,w)}:~k\geq 1\}_{z,w\not\in B, z\neq w}, \{\hat{\eta}_0^N(z)\}_{z\in B}, \{\eta_0^N(z)\}_{z\not\in B}\Bigg)
\end{align*}
is independent of $1_{\{\Gamma_{t,x}=B\}}$ and
\[
Q_{x, B}\left(\{\Xi_s^{(z,w)}:~0\leq s\leq t\}_{(z,w)\in \beta_B}, \{U^{(z,w)}_k:~k\geq 1\}_{(z,w)\in \beta_B}, \{\eta_0^N(z)\}_{z\in B} \right).
\]
As a result,
\[
\mathbb{P}\left(\widetilde{\eta}_t^N(y)=k_1, \eta_t^N(x)=k_2\Big|\Gamma_{t,x}=B\right)=A(k_1, B)\mathbb{P}\left(\eta_t^N(x)=k_2\Big|\Gamma_{t,x}=B\right),
\]
where $A(k_1, B)$  is the probability of the event
\begin{align*}
&\widetilde{Q}_y\Bigg(\{\hat{\Xi}_t^{(z,w)}:~t\geq 0\}_{z\in B\text{~or~}w\in B, z\neq w}, \{\Xi_t^{(z,w)}:~t\geq 0\}_{z,w\not\in B, z\neq w},  \\
&\text{\quad}\{\hat{U}_k^{(z,w)}:~k\geq 1\}_{z\in B\text{~or~}w\in B, z\neq w}, \{U_k^{(z,w)}:~k\geq 1\}_{z,w\not\in B, z\neq w}, \{\hat{\eta}_0^N(z)\}_{z\in B}, \{\eta_0^N(z)\}_{z\not\in B}\Bigg)\\
&=k_1.
\end{align*}
Since $\{\{\hat{\Xi}_t^{(z,w)}:~t\geq 0\}_{z\neq w}, \{\hat{U}_k^{(z,w)}:~k\geq 1\}_{z\neq w}, \{\hat{\eta}_0^N(z)\}_{1\leq z\leq N}\}$ is an independent copy of $\{\{\Xi_t^{(z,w)}:~t\geq 0\}_{z\neq w}, \{U_k^{(z,w)}:~k\geq 1\}_{z\neq w}, \{\eta_0^N(z)\}_{1\leq z\leq N}\}$ and all these Poisson processes, uniformly distributed random variables, $\{\eta_0^N(z)\}_{1\leq z\leq N}$, $\{\hat{\eta}_0^N(z)\}_{1\leq z\leq N}$ are independent, we have that $A(k_1, B)$ is independent of $B$. Consequently, Equation \eqref{equ 3.1} holds and the proof is complete.

\qed

\proof[Proof of Lemma \ref{lemma 3.1.3}]

We first bound from above the probability that $w\in \Gamma_{t,z}$ for any $w\neq z$. For integer $m\geq 1$ and a given path $x_0, x_1,\ldots,x_m$ with $x_0=w$ and $x_m=z$, the probability that $x_0, x_1,\ldots,x_m$ is a $t$-influence path equals
\[
\mathbb{P}(\sum_{i=1}^m\tau_i\leq t),
\]
where $\tau_1,\ldots,\tau_m$ are independent and identically distributed with the exponential distribution with rate $2K_1/N$. Using Markov inequality,
\[
\mathbb{P}(\sum_{i=1}^m\tau_i\leq t)\leq e^{\theta t}\left(\frac{2K_1}{N\theta+2K_1}\right)^m
\]
for any $\theta>0$. Let $\theta=4K_1$, then $\mathbb{P}(\sum_{i=1}^m\tau_i\leq t)\leq e^{4K_1t}\frac{1}{2^m}\frac{1}{N^m}$. The number of paths from $w$ to $z$ with length $m$ is at most $N^{m-1}$, hence
\begin{align}\label{equ 3.2}
\mathbb{P}(w\in \Gamma_{t,z})\leq \sum_{m=1}^{+\infty}e^{4K_1t}\frac{1}{2^m}\frac{N^{m-1}}{N^m}=\frac{1}{N}e^{4K_1t}\leq \frac{1}{N}e^{4K_1T}
\end{align}
for $0\leq t\leq T$. For $z\neq x,y$, we denote by $\{z\in \Gamma_{t,x}\}\circ\{z\in \Gamma_{t,y}\}$ the event that there are two $t$-influence paths from $z$ to $x$ and $y$ respectively such that these two paths do not have common edges. According to Berg-Kesten inequality and Equation \eqref{equ 3.2},
\begin{equation}\label{equ 3.3}
\mathbb{P}(\{z\in \Gamma_{t,x}\}\circ\{z\in \Gamma_{t,y}\})\leq
\mathbb{P}(z\in \Gamma_{t,x})\mathbb{P}(z\in \Gamma_{t,y})\leq \frac{e^{8K_1T}}{N^2}.
\end{equation}
If $x\not\in \Gamma_{t,y}$ and $y\not\in \Gamma_{t,x}$, then $\Gamma_{t,x}\bigcap \Gamma_{t,y}\neq \emptyset$ implies that there exists $z\neq x,y$ such that $\{z\in \Gamma_{t,x}\}\circ\{z\in \Gamma_{t,y}\}$ occurs. Hence, by Equations \eqref{equ 3.2} and \eqref{equ 3.3},
\begin{align*}
\mathbb{P}(\Gamma_{t,x}\bigcap \Gamma_{t,y}\neq \emptyset)&\leq \mathbb{P}(x\in \Gamma_{t,y})+\mathbb{P}(y\in \Gamma_{t,x})
+\sum_{z\neq x,y}\mathbb{P}(\{z\in \Gamma_{t,x}\}\circ\{z\in \Gamma_{t,y}\})\\
&\leq \frac{2}{N}e^{4K_1T}+\frac{(N-2)e^{8K_1T}}{N^2}.
\end{align*}
Consequently, Lemma \ref{lemma 3.1.3} holds with $C_3=2e^{4K_1T}+e^{8K_1T}$.

\qed

At last, we prove Lemma \ref{lemma 3.1 approximated independence finite K}.

\proof[Proof of Lemma \ref{lemma 3.1 approximated independence finite K}]

For any given $z\in \{1,2,\ldots, N\}$ and $\{1,\ldots,N\}\supseteq C\not\ni z$, we denote by $\Gamma_{t,z}^{\setminus C}$ the set
\[
\{z\}\bigcup\{w:~\text{there is a $t$-influence path from $w$ to $z$ avoiding $C$}\}.
\]
We denote by $D_t$ the event that
\[
\hat{\Xi}_t^{(w_2, w_1)}=\hat{\Xi}_t^{(w_1,w_2)}=0
\]
for all $w_1\in \Gamma_{t,x}^{\setminus \{y\}}$ and $w_2\in \Gamma_{t,y}^{\setminus(\Gamma_{t,x}^{\setminus\{y\}})}$. For given $z\neq w$,
\[
\mathbb{P}(\hat{\Xi}_t^{(z, w)}>0\text{~or~}\hat{\Xi}_t^{(w,z)}>0)\leq 2(1-e^{-K_1t/N})\leq \frac{2K_1T}{N}.
\]
Since $\{\hat{\Xi}_t^{(z,w)}\}_{z\neq w}$ and $\{\Xi_t^{(z,w)}\}_{z\neq w}$ are independent, we have
\begin{equation}\label{equ 3.4}
\mathbb{P}(D_t^c)\leq \frac{2K_1T}{N} \mathbb{E}\left(\left|\Gamma_{t,x}^{\setminus \{y\}}\right|\left|\Gamma_{t,y}^{\setminus(\Gamma_{t,x}^{\setminus\{y\}})}\right|\right),
\end{equation}
where $|A|$ is the cardinality of the set $A$. On the event $D_t\bigcap\{\Gamma_{t,x}\bigcap\Gamma_{t,y}=\emptyset\}$, we have
\[
\Gamma_{t,x}=\Gamma_{t,x}^{\setminus \{y\}}, \text{\quad}
\Gamma_{t,y}=\Gamma_{t,y}^{\setminus(\Gamma_{t,x}^{\setminus\{y\}})}=\widetilde{\Gamma}_{t,y}
\]
and then $\eta_t^N(y)=\widetilde{\eta}_t^N(y)$.  Hence, for $0\leq k_1, k_2\leq \mathcal{K}$,
\begin{equation}\label{equ 3.5}
\left|\mathbb{P}(\eta_t^N(y)=k_2)-\mathbb{P}(\widetilde{\eta}_t^N(y)=k_2)\right|\leq 2\left(\mathbb{P}(D_t^c)+\mathbb{P}(\Gamma_{t,x}\bigcap\Gamma_{t,y}\neq\emptyset)\right)
\end{equation}
and
\begin{align*}
&\left|\mathbb{P}(\eta_t^N(x)=k_1, \eta_t^N(y)=k_2)-\mathbb{P}(\eta_t^N(x)=k_1, \widetilde{\eta}_t^N(y)=k_2)\right| \\
&\leq 2\left(\mathbb{P}(D_t^c)+\mathbb{P}(\Gamma_{t,x}\bigcap\Gamma_{t,y}\neq\emptyset)\right).
\end{align*}
Then, by Lemma \ref{lemma 3.1.2} and Equation \eqref{equ 3.5},
\begin{align}\label{equ 3.55}
&\left|\mathbb{P}\left(\eta_t^N(x)=k_1, \eta_t^N(y)=k_2\right)-\mathbb{P}(\eta_t^N(x)=k_1)\mathbb{P}(\eta_t^N(y)=k_2)\right| \notag\\
&\leq \left|\mathbb{P}\left(\eta_t^N(x)=k_1, \eta_t^N(y)=k_2\right)-\mathbb{P}(\eta_t^N(x)=k_1, \widetilde{\eta}_t^N(y)=k_2)\right| \notag\\
&\text{\quad}+\left|\mathbb{P}(\eta_t^N(x)=k_1)\mathbb{P}(\widetilde{\eta}_t^N(y)=k_2)-\mathbb{P}(\eta_t^N(x)=k_1)\mathbb{P}(\eta_t^N(y)=k_2)\right| \notag\\
&\leq 4\left(\mathbb{P}(D_t^c)+\mathbb{P}(\Gamma_{t,x}\bigcap\Gamma_{t,y}\neq\emptyset)\right).
\end{align}
Using the total probability formula,
\begin{align*}
\mathbb{E}\left(\left|\Gamma_{t,x}^{\setminus \{y\}}\right|\left|\Gamma_{t,y}^{\setminus(\Gamma_{t,x}^{\setminus\{y\}})}\right|\right)=\sum_{B\not\ni y, B\ni x}\mathbb{E}\left(|\Gamma_{t,y}^{\setminus B}|\Big|\Gamma_{t,x}^{\setminus \{y\}}=B\right)|B|\mathbb{P}(\Gamma_{t,x}^{\setminus \{y\}}=B).
\end{align*}
According to an analysis similar with that given in the proof of Lemma \ref{lemma 3.1.1}, $\{\Gamma_{t,x}^{\setminus \{y\}}=B\}$ is independent of $\{\Xi_s^{(z,w)}:0\leq s\leq t\}_{z\neq w, z,w\not\in B}$. Hence,
\[
\mathbb{E}\left(|\Gamma_{t,y}^{\setminus B}|\Big|\Gamma_{t,x}^{\setminus \{y\}}=B\right)
=\mathbb{E}\left(|\Gamma_{t,y}^{\setminus B}|\right)\leq \mathbb{E}\left(|\Gamma_{t,y}|\right)
\]
and consequently,
\begin{align}\label{equ 3.6}
&\mathbb{E}\left(\left|\Gamma_{t,x}^{\setminus \{y\}}\right|\left|\Gamma_{t,y}^{\setminus(\Gamma_{t,x}^{\setminus\{y\}})}\right|\right) \leq \mathbb{E}\left(|\Gamma_{t,y}|\right)\sum_{B\not\ni y, B\ni x}|B|\mathbb{P}(\Gamma_{t,x}^{\setminus \{y\}}=B)\notag\\
&=\mathbb{E}\left(|\Gamma_{t,y}|\right)
\mathbb{E}\left(|\Gamma_{t,x}^{\setminus \{y\}}|\right)\leq \mathbb{E}\left(|\Gamma_{t,x}|\right)\mathbb{E}\left(|\Gamma_{t,y}|\right).
\end{align}

Using Equation \eqref{equ 3.2}, for any $z\in \{1,2,\ldots, N\}$,
\[
\mathbb{E}\left(|\Gamma_{t,z}|\right)\leq 1+\frac{N-1}{N}e^{4K_1T}\leq 1+e^{4K_1T}.
\]
Then, by Equations \eqref{equ 3.4} and \eqref{equ 3.6},
\[
\mathbb{P}(D_t^c)\leq \frac{2K_1T}{N}\left(1+e^{4K_1T}\right)^2.
\]
Finally, using Equation \eqref{equ 3.55} and Lemma \ref{lemma 3.1.3}, Lemma \ref{lemma 3.1 approximated independence finite K} holds with
\[
C_2=4\left(2K_1T\left(1+e^{4K_1T}\right)^2+C_3\right),
\]
where $C_3$ is defined as in Lemma \ref{lemma 3.1.3}.

\qed

\subsection{The proof of Lemma \ref{lemma 3.2 approximated independence infinite K}}\label{subsection 3.2}

In this subsection we prove Lemma \ref{lemma 3.2 approximated independence infinite K}. Throughout this subsection, we assume that $T>0$ is a given moment. We first introduce some definitions and notations for later use. For $N\geq 1$, we denote by $\{\zeta_t^N\}_{t\geq 0}$ the Markov process with state space $\{0,1,2,\ldots\}^N$ and generator $\Omega_N$ given by
\begin{equation}\label{equ 3.11 generator of auxiliary process}
\Omega_Nf(\eta)=\frac{1}{N}\sum_{i=1}^N\sum_{j\neq i}\sum_{k\in \mathbb{Z}_\mathcal{K}}\sum_{l\in \mathbb{Z}_{\mathcal{K}}}
\hat{\phi}_{k,l}\left(i/N, j/N\right)1_{\{\eta(i)=k, \eta(j)=l\}}\left(f(\eta^{(i,j)})-f(\eta)\right)
\end{equation}
for any $\eta\in \{0,1,\ldots\}^N$ and $f$ from $\{0,1,\ldots\}^N$ to $\mathbb{R}$, where
\[
\hat{\phi}_{k,l}=
\begin{cases}
\phi_{k,l} &\text{~if~}k\leq \lfloor C_4\log N\rfloor,\\
\phi_{\lfloor C_4\log N\rfloor, l} &\text{~if~} k>\lfloor C_4\log N\rfloor,
\end{cases}
\]
where $C_4=\frac{1}{250C_1T}$ and $C_1$ is defined as in Assumption (B). We further assume that $\zeta_0^N=\eta_0^N$.

The following lemma is crucial for us to prove Lemma \ref{lemma 3.2 approximated independence infinite K}.

\begin{lemma}\label{lemma 3.2.1}
Under Assumptions (A) and (B), there exists $N_1\geq 1$ such that
\[
\left|\mathbb{P}(\eta_t^N(x)=k_1)-\mathbb{P}(\zeta_t^N(x)=k_1)\right|\leq \frac{1}{N^{0.95}}
\]
and
\[
\left|\mathbb{P}(\eta_t^N(x)=k_1, \eta_t^N(y)=k_2)-\mathbb{P}(\zeta_t^N(x)=k_1, \zeta_t^N(y)=k_2)\right|\leq \frac{1}{N^{0.95}}
\]
for any $N\geq N_1, 0\leq t\leq T, k_1,k_2\geq 0$ and any two different $x,y\in \{1,2,\ldots,N\}$.
\end{lemma}

We prove Lemma \ref{lemma 3.2.1} at the end of this subsection. Now we utilize Lemma \ref{lemma 3.2.1} to prove Lemma \ref{lemma 3.2 approximated independence infinite K}.

\proof[Proof of Lemma \ref{lemma 3.2 approximated independence infinite K}]

According to a graphical analysis similar with that given in the proof of Lemma \ref{lemma 3.1 approximated independence finite K}, we have
\begin{equation}\label{equ 3.12}
\left|\mathbb{P}(\zeta_t(x)=k_1, \zeta_t(y)=k_2)-\mathbb{P}(\zeta_t(x)=k_1)\mathbb{P}(\zeta_t(y)=k_2)\right|\leq \frac{\hat{C}_2}{N}
\end{equation}
for any $N\geq 1$, $k_1,k_2\geq 0, 0\leq t\leq T$ and two different $x,y\in \{1,\ldots,N\}$, where
\[
\hat{C}_2=4\left(2\hat{K}_1T(1+e^{4\hat{K}_1T})^2+\hat{C}_3\right), \text{\quad} \hat{K}_1=\sup_{k,l\geq 0, u,v\in T}\hat{\phi}_{k,l}(u,v)
\]
and
\[
\hat{C}_3=2e^{4\hat{K}_1T}+e^{8\hat{K}_1T}.
\]
By Assumption (B) and the definition of $\{\hat{\phi}_{k,l}\}_{k,l\geq 0}$, we have $\hat{K}_1\leq C_1C_4\log N$. Then, by the definition of $C_4$, there exists
$N_2\geq 1$ such that
\begin{equation*}
\frac{\hat{C}_2}{N}\leq \frac{1}{N^{0.95}}
\end{equation*}
when $N\geq N_2$. Therefore, by Lemma \ref{lemma 3.2.1} and Equation \eqref{equ 3.12},
\[
\left|\mathbb{P}\left(\eta_t^N(x)=k_1, \eta_t^N(y)=k_2\right)-\mathbb{P}(\eta_t^N(x)=k_1)\mathbb{P}(\eta_t^N(y)=k_2)\right|\leq \frac{\hat{C}_2}{N}+\frac{3}{N^{0.95}}\leq \frac{4}{N^{0.95}}
\]
for any $N\geq \max\{N_1, N_2\}$, $k_1,k_2\geq 0, 0\leq t\leq T$ and two different $x,y\in \{1,\ldots,N\}$. Since $\frac{4}{N^{0.95}}\leq \frac{1}{N^{0.9}}$ for sufficiently large $N$, the proof is complete.

\qed

At last, we prove Lemma \ref{lemma 3.2.1}.

\proof[Proof of Lemma \ref{lemma 3.2.1}]
Let
\[
\varpi^N=\inf\{t:~\eta_t^N(x)\geq C_4\log N\text{~for some~}x\in\{1,\ldots,N\}\}.
\]
Since $\zeta_0^N=\eta_0^N$ and $\phi_{k,l}\neq \hat{\phi}_{k,l}$ only when $k\geq C_4\log N$, we have
\[
\eta_t^N=\zeta_t^N
\]
for $0\leq t\leq \varpi^N$ in the sense of coupling.  Under Assumption (A), using Markov inequality,
\[
\mathbb{P}\left(\eta_0^N(x)>\frac{C_4\log N}{5}\right)\leq e^{-\theta_1\frac{C_4\log N}{5}}e^{\psi(x/N)(e^{\theta_1}-1)}
\leq N^{-\frac{C_4\theta_1}{5}}e^{\|\psi\|_\infty(e^{\theta_1}-1)}
\]
and
\begin{align*}
\mathbb{P}\left(\sum_y\eta_0^N(y)>3N\int_\mathbb{T}\psi(u)du\right)&\leq e^{-3\theta_2N\int_\mathbb{T}\psi(u)du}\prod_ye^{(e^{\theta_2}-1)\psi(y/N)}\\
&=e^{-N\int_\mathbb{T}\psi(u)du(3\theta_2-(e^{\theta_2}-1)+o(1))}
\end{align*}
for any $\theta_1, \theta_2>0$. Let $\theta_1=\frac{15}{C_4}$ and $\theta_2=\log 3$, then
\begin{equation}\label{equ 3.13}
\mathbb{P}\left(\eta_0^N(x)>\frac{C_4\log N}{5}\right)\leq N^{-3}e^{\|\psi\|_\infty(e^{15/C_4}-1)}
\end{equation}
and
\begin{equation}\label{equ 3.14}
\mathbb{P}\left(\sum_y\eta_0^N(y)>3N\int_\mathbb{T}\psi(u)du\right)
\leq e^{-N\int_\mathbb{T}\psi(u)du(\log 27-2+o(1))}.
\end{equation}

Conditioned on $\sum_y\eta_0^N(y)\leq 3N\int_\mathbb{T}\psi(u)du$, using Assumption (B), $\sup_{0\leq t\leq T}(\eta_t^N(x)-\eta_0^N(x))$ is stochastically dominated from above by $\mathcal{N}_T$, where $\{\mathcal{N}_t\}_{t\geq 0}$ is a Poisson process with rate
\[
3N\int_\mathbb{T}\psi(u)du\frac{C_1}{N}=3C_1\int_\mathbb{T}\psi(u)du.
\]
Hence, using Markov inequality,
\begin{align*}
&\mathbb{P}\left(\sup_{0\leq t\leq T}(\eta_t^N(x)-\eta_0^N(x))\geq \frac{C_4\log N}{5}\Bigg|\sum_y\eta_0^N(y)\leq 3N\int_\mathbb{T}\psi(u)du\right)\\
&\leq \mathbb{P}\left(\mathcal{N}_T\geq \frac{C_4\log N}{5}\right)\leq e^{-\frac{C_4\theta_3\log N}{5}}e^{3C_1\int_\mathbb{T}\psi(u)du(e^{\theta_3}-1)}\\
&=N^{-\frac{C_4\theta_3}{5}}e^{3C_1(e^{\theta_3}-1)\int_\mathbb{T}\psi(u)du}.
\end{align*}
Let $\theta_3=\frac{15}{C_4}$, then
\begin{align*}
&\mathbb{P}\left(\sup_{0\leq t\leq T}(\eta_t^N(x)-\eta_0^N(x))\geq \frac{C_4\log N}{5}\Bigg|\sum_y\eta_0^N(y)\leq 3N\int_\mathbb{T}\psi(u)du\right)\\
&\leq N^{-3}e^{3C_1(e^{\frac{15}{C_4}}-1)\int_\mathbb{T}\psi(u)du}.
\end{align*}
As a result, by Equation \eqref{equ 3.14},
\begin{align*}
&\mathbb{P}\left(\sup_{0\leq t\leq T}(\eta_t^N(x)-\eta_0^N(x))\geq \frac{C_4\log N}{5}\right)\\
&\leq  N^{-3}e^{3C_1(e^{\frac{15}{C_4}}-1)\int_\mathbb{T}\psi(u)du}
+e^{-N\int_\mathbb{T}\psi(u)du(\log 27-2+o(1))}\leq N^{-2}
\end{align*}
for sufficiently large $N$. Consequently, by Equation \eqref{equ 3.13},
\[
\mathbb{P}\left(\sup_{0\leq t\leq T}\eta_t^N(x)\geq \frac{2C_4}{5}\log N\right)
\leq N^{-3}e^{\|\psi\|_\infty(e^{15/C_4}-1)}+N^{-2}\leq 2N^{-2}
\]
for sufficiently large $N$. Therefore,
\begin{equation}\label{equ 3.15}
\mathbb{P}(\varpi^N\leq T)\leq 2N^{-2}N=2N^{-1}\leq \frac{1}{N^{0.99}}
\end{equation}
for sufficiently large $N$. As we have explained, $\eta_t^N=\zeta_t^N$ for $0\leq t\leq \varpi^N$. Therefore,
\[
\left|\mathbb{P}(\eta_t^N(x)=k_1)-\mathbb{P}(\zeta_t^N(x)=k_1)\right|\leq 2\mathbb{P}(\varpi^N\leq T)
\]
and
\[
\left|\mathbb{P}(\eta_t^N(x)=k_1, \eta_t^N(y)=k_2)-\mathbb{P}(\zeta_t^N(x)=k_1, \zeta_t^N(y)=k_2)\right|\leq 2\mathbb{P}(\varpi^N\leq T)
\]
for any $0\leq t\leq T, k_1,k_2\geq 0$ and any two different $x,y\in \{1,2,\ldots,N\}$. Then, Lemma \ref{lemma 3.2.1} follows from Equation \eqref{equ 3.15} and the proof is complete.

\qed

\section{The proof of Theorem \ref{theorem 2.1 mean field limit finite}}\label{section four proof of hydro finite}
In this section we prove Theorem \ref{theorem 2.1 mean field limit finite}. Throughout this section, we assume $\mathcal{K}<+\infty$ and adopt Assumption (A). For later use, we introduce some notations and definitions. For any $f=(f_0,f_1,\ldots,f_\mathcal{K})\in \left(C(\mathbb{T})\right)^{\mathcal{K}+1}$, we define
\[
\|f\|_\infty^{\mathcal{K}}=\sup_{0\leq k\leq \mathcal{K}, u\in \mathbb{T}}|f_k(u)|.
\]
For any $t\geq 0$, $0\leq k\leq \mathcal{K}$ and $u\in \mathbb{T}$, we define
\[
\varrho_{t,k}^N(u)=\mathbb{P}\left(\eta_t^N(i)=k\right)
\]
when $\frac{i-1}{N}<u\leq \frac{i}{N}$ for some $i\in \{1,2,\ldots, N\}$.

As preliminaries, we first give two lemmas.
\begin{lemma}\label{lemma 4.1}
Let $\left(C(\mathbb{T})\right)^{\mathcal{K}+1}$ be endowed with the norm $\|\cdot\|_\infty^\mathcal{K}$, then the solution
\[
\{\rho_k^{\mathcal{K}}(t,u):~t\geq 0, u\in \mathbb{T}\}_{0\leq k\leq \mathcal{K}}
\]
to Equation \eqref{equ 2.1 mean field ODE finite types} exists on an unbounded area of $[0, +\infty)\times \left(C(\mathbb{T})\right)^{\mathcal{K}+1}$ and is unique.
\end{lemma}

\proof

Under the norm $\|\cdot\|_\infty^\mathcal{K}$, Equation \eqref{equ 2.1 mean field ODE finite types} satisfies the local Lipschitz condition. Hence, Lemma \ref{lemma 4.1} follows from Theorem 19.1.2 of \cite{Lang}.

\qed

Note that Lemma \ref{lemma 4.1} does not ensure that the solution to Equation \eqref{equ 2.1 mean field ODE finite types} exists for all $t\geq 0$, since this lemma does not exclude the case where the solution exists for $t<T_0$ for some $T_0$ and $\lim_{t\uparrow T_0}\|\{\rho_{t,k}^\mathcal{K}\}_{0\leq k\leq \mathcal{K}}\|_\infty^\mathcal{K}=+\infty$.

\begin{lemma}\label{lemma 4.2}
For any $T_1>0$, if the solution $\{\rho_{t,k}^{\mathcal{K}}:~t\geq 0\}_{0\leq k\leq \mathcal{K}}$ to Equation \eqref{equ 2.1 mean field ODE finite types} exists for $0\leq t\leq T_1$, then
\[
\lim_{N\rightarrow+\infty}\sup_{0\leq k\leq \mathcal{K}, u\in \mathbb{T}}\left|\rho_{t,k}^\mathcal{K}(u)-\varrho_{t,k}^N(u)\right|=0
\]
for all $0\leq t\leq T_1$.
\end{lemma}

\proof

Utilizing the generator $\mathcal{L}_N$ of $\{\eta_t^N\}_{t\geq 0}$ given in Equation \eqref{equ 1.1 generator} and Kolmogorov-Chapman equation, for $0\leq t\leq T_1$, $1\leq k\leq \mathcal{K}-1$ and $u\in \big(\frac{i-1}{N}, \frac{i}{N}\big]$ for some $i\in \{1,\ldots,N\}$, we have
\begin{align*}
\frac{d}{dt}\varrho_{t,k}^N(u)=&-\frac{1}{N}\sum_{l=0}^\mathcal{K}\sum_{j=1}^N\phi_{k,l}(i/N, j/N)\mathbb{P}\left(\eta_t^N(i)=k, \eta_t^N(j)=l\right)\\
&-\frac{1}{N}\sum_{l=1}^\mathcal{K}\sum_{j=1}^N\phi_{l,k}(j/N,i/N)\mathbb{P}\left(\eta_t^N(i)=k, \eta_t^N(j)=l\right)\\
&+\frac{1}{N}\sum_{l=1}^\mathcal{K}\sum_{j=1}^N\phi_{l,k-1}(j/N,i/N)\mathbb{P}\left(\eta_t^N(i)=k-1, \eta_t^N(j)=l\right)\\
&+\frac{1}{N}\sum_{l=0}^\mathcal{K}\sum_{j=1}^N\phi_{k+1,l}(i/N,j/N)\mathbb{P}\left(\eta_t^N(i)=k+1, \eta_t^N(j)=l\right).
\end{align*}
By Lemma \ref{lemma 3.1 approximated independence finite K}, when we replace $\mathbb{P}\left(\eta_t^N(i)=k_1, \eta_t^N(j)=k_2\right)$ by
\[
\mathbb{P}\left(\eta_t^N(i)=k_1\right)\mathbb{P}\left(\eta_t^N(j)=k_2\right)
\]
for $i\neq j$, the error is at most $\frac{C_2}{N}$. Hence,
\begin{align*}
\frac{d}{dt}\varrho_{t,k}^N(u)=&-\frac{1}{N}\varrho_{t,k}^N(u)\sum_{l=0}^\mathcal{K}\sum_{j=1}^N\phi_{k,l}(i/N, j/N)\varrho_{t,l}^N(j/N)\\
&-\frac{1}{N}\varrho_{t,k}^N(u)\frac{1}{N}\sum_{l=1}^\mathcal{K}\sum_{j=1}^N\phi_{l,k}(j/N,i/N)\varrho_{t,l}^N(j/N)\\
&+\varrho_{t, k-1}^N(u)\frac{1}{N}\sum_{l=1}^\mathcal{K}\sum_{j=1}^N\phi_{l,k-1}(j/N,i/N)\varrho_{t,l}^N(j/N)\\
&+\varrho_{t,k+1}^N(u)\frac{1}{N}\sum_{l=0}^\mathcal{K}\sum_{j=1}^N\phi_{k+1,l}(i/N,j/N)\varrho_{t,l}^N(j/N)+\varepsilon_{t,1}^N,
\end{align*}
where
\begin{equation*}
\sup_{0\leq t\leq T_1}|\varepsilon_{t,1}^N|=O(N^{-1}).
\end{equation*}
Since $\phi_{k,l}\in C^\infty(\mathbb{T}^2)$, we have
\[
\sup\left\{\left|\phi_{k,l}(u_1, v_1)-\phi_{k,l}(u_2, v_2)\right|:~|u_1-u_2|\leq N^{-1}, |v_1-v_2|\leq N^{-1}\right\}=O(N^{-1}).
\]
Hence,
\begin{align}\label{equ 4.2}
\frac{d}{dt}\varrho_{t,k}^N(u)=&-\varrho_{t,k}^N(u)\sum_{l=0}^\mathcal{K}\int_{v\in \mathbb{T}}\phi_{k,l}(u, v)\varrho_{t,l}^N(v)dv\notag\\
&-\varrho_{t,k}^N(u)\sum_{l=1}^\mathcal{K}\int_{v\in \mathbb{T}}\phi_{l,k}(v, u)\varrho_{t,l}^N(v)dv\\
&+\varrho_{t, k-1}^N(u)\sum_{l=1}^\mathcal{K}\int_{v\in \mathbb{T}}\phi_{l,k-1}(v,u)\varrho_{t,l}^N(v)dv\notag\\
&+\varrho_{t,k+1}^N(u)\sum_{l=0}^\mathcal{K}\int_{v\in \mathbb{T}}\phi_{k+1,l}(u,v)\varrho_{t,l}^N(v)dv+\varepsilon_{t,2}^N,\notag
\end{align}
where
\begin{equation*}
\sup_{0\leq t\leq T_1}|\varepsilon_{t,2}^N|=O(N^{-1}).
\end{equation*}
According to similar analyses, we have
\begin{align}\label{equ 4.3}
\frac{d}{dt}\varrho_{t,0}^N(u)=&-\varrho_{t,0}^N(u)\sum_{l=1}^\mathcal{K}\int_{\mathbb{T}}\phi_{l,0}(v,u)\varrho_{t,l}^N(v)dv\notag\\
&+\varrho_{t,1}^N(u)\sum_{l=0}^{\mathcal{K}}\int_{\mathbb{T}}\phi_{1,l}(u,v)\varrho_{t,l}^N(v)dv+\varepsilon_{t,3}^N
\end{align}
and
\begin{align}\label{equ 4.4}
\frac{d}{dt}\varrho_{t,\mathcal{K}}^N(u)=&-\varrho_{t,\mathcal{K}}^N(u)\sum_{l=0}^\mathcal{K}\int_{\mathbb{T}}\phi_{\mathcal{K},l}(u,v)\varrho_{t,l}^N(v)dv\notag\\
&+\varrho_{t,\mathcal{K}-1}^N(u)\sum_{l=1}^{\mathcal{K}}\int_{\mathbb{T}}\phi_{l,\mathcal{K}-1}(v,u)\varrho_{t,l}^N(v)dv+\varepsilon_{t,4}^N,
\end{align}
where
\begin{equation*}
\sup_{0\leq t\leq T_1}|\varepsilon_{t,3}^N|=O(N^{-1})\text{~and~}\sup_{0\leq t\leq T_1}|\varepsilon_{t,4}^N|=O(N^{-1}).
\end{equation*}
According to Assumption (A), we have
\begin{equation}\label{equ 4.5}
\varepsilon_5^N:=\sup_{0\leq k\leq \mathcal{K}, u\in \mathbb{T}}\left|\rho_{0,k}^\mathcal{K}(u)-\varrho_{0,k}^N(u)\right|=O(N^{-1}).
\end{equation}
Let $C_5=\sup_{0\leq t\leq T_1, 0\leq k\leq \mathcal{K}, u\in \mathbb{T}}|\rho_{t,k}^{\mathcal{K}}(u)|$.
Using triangle inequality, Grownwall's inequality and Equations \eqref{equ 4.2}-\eqref{equ 4.5}, we have
\[
\sup_{0\leq k\leq \mathcal{K}, u\in \mathbb{T}}\left|\rho_{t,k}^\mathcal{K}(u)-\varrho_{t,k}^N(u)\right|\leq \varepsilon_6^N e^{C_6t}
\]
for $0\leq t\leq T_1$, where
\[
C_6=4(\mathcal{K}+1)\sup_{0\leq k,l\leq \mathcal{K}}\|\phi_{k,l}\|_\infty(C_5+1)
\]
and
\[
\varepsilon_6^N=\varepsilon_5^N+T_1\left(\sum_{l=2}^4\sup_{0\leq t\leq T_1}|\varepsilon_{t,l}^N|\right)=O(N^{-1}).
\]
Hence, Lemma \ref{lemma 4.2} holds and the proof is complete.

\qed

At last, we prove Theorem \ref{theorem 2.1 mean field limit finite}.

\proof[Proof of Theorem \ref{theorem 2.1 mean field limit finite}]

By Lemma \ref{lemma 4.2}, for any $T_1>0$, if the solution $\{\rho_{t,k}^{\mathcal{K}}:~t\geq 0\}_{0\leq k\leq \mathcal{K}}$ to Equation \eqref{equ 2.1 mean field ODE finite types} exists for $0\leq t\leq T_1$, then $0\leq \rho_{t,k}^{\mathcal{K}}(u)\leq 1$ for any $0\leq k\leq \mathcal{K}, u\in \mathbb{T}$ and $0\leq t\leq T_1$. Consequently, using Lemma \ref{lemma 4.1}, the solution $\{\rho_{t,k}^{\mathcal{K}}:~t\geq 0\}_{0\leq k\leq \mathcal{K}}$ to Equation \eqref{equ 2.1 mean field ODE finite types} exists for all $t\geq 0$ and is unique. Then, by Lemma \ref{lemma 4.2},
\[
\lim_{N\rightarrow+\infty}\mathbb{E}\mu_{t,k}^N(f)=\int_\mathbb{T}\rho_{t,k}^\mathcal{K}(u)f(u)du
\]
for any $t\geq 0$, $0\leq k\leq \mathcal{K}$ and $f\in C(\mathbb{T})$. Hence, to prove Theorem \ref{theorem 2.1 mean field limit finite} we only need to show that
\begin{equation}\label{equ 4.1}
\lim_{N\rightarrow+\infty}{\rm Var}\left(\mu_{t,k}^N(f)\right)=0
\end{equation}
for any $t\geq 0$, $0\leq k\leq \mathcal{K}$ and $f\in C(\mathbb{T})$. Using Lemma \ref{lemma 3.1 approximated independence finite K}, we have
\begin{align*}
{\rm Var}\left(\mu_{t,k}^N(f)\right)
&=\frac{1}{N^2}\sum_{x=1}^N\sum_{y=1}^N{\rm Cov}\left(1_{\{\eta_t^N(x)=k\}}, 1_{\{\eta_t^N(y)=k\}}\right)f(x/N)f(y/N)\\
&\leq \frac{2N}{N^2}\|f\|_\infty^2+\frac{1}{N^2}\sum_{x=1}^N\sum_{y\neq x}{\rm Cov}\left(1_{\{\eta_t^N(x)=k\}}, 1_{\{\eta_t^N(y)=k\}}\right)f(x/N)f(y/N)\\
&\leq \frac{2N}{N^2}\|f\|_\infty^2+\frac{N(N-1)}{N^2}\frac{C_2}{N}\|f\|_\infty^2=O(N^{-1}).
\end{align*}
As a result, Equation \eqref{equ 4.1} holds and the proof is complete.

\qed

\section{The proof of Theorem \ref{theorem 2.2 mean field limit infinite}}\label{section five proof of hydro infinite}
In this section, we prove Theorem \ref{theorem 2.2 mean field limit infinite}. As preliminaries, we first give several lemmas.

\begin{lemma}\label{lemma 5.1}
Under Assumptions (A) and (B),
\begin{equation}\label{equ 5.1}
\sup_{0\leq t\leq T, 1\leq x\leq N}\mathbb{E}\eta_t^N(x)\leq (TC_1+1)\|\psi\|_\infty,
\end{equation}
\begin{align}\label{equ 5.1 two}
\sup_{0\leq t\leq T, 1\leq x\leq N}\mathbb{E}\left((\eta_t^N(x))^2\right)
\leq 2(\|\psi\|_\infty^2+\|\psi\|_\infty)+2TC_1\|\psi\|_\infty+2T^2C_1^2\left(\|\psi\|^2_\infty+\|\psi\|_\infty\right)
\end{align}
for all $N\geq 1$ and
\begin{equation}\label{equ 5.2}
\limsup_{M\rightarrow+\infty}\frac{1}{M}\sup_{0\leq t\leq T, N\geq 1, 1\leq x\leq N}\log\left(\sum_{l=M}^{+\infty}l\mathbb{P}(\eta_t^N(x)=l)\right)=-\infty.
\end{equation}
\end{lemma}

\proof

Conditioned on $\eta_0^N$, under Assumption (B)$, \eta_t^N(x)$ is stochastically dominated from above by
\[
\eta_0^N(x)+\mathcal{Y}\left(t\frac{C_1}{N}\sum_{y=1}^N\eta_0^N(y)\right),
\]
where $\{\mathcal{Y}_t\}_{t\geq 0}$ is a Poisson process with rate $1$. Hence, under Assumption (A),
\[
\mathbb{E}\eta_t^N(x)\leq \mathbb{E}\eta_0^N(x)+t\frac{C_1}{N}\sum_{y=1}^N\mathbb{E}\eta_0^N(y)\leq \|\psi\|_\infty+tC_1\|\psi\|_\infty
\]
and consequently Equation \eqref{equ 5.1} holds. Similarly,
\[
\mathbb{E}\left((\eta_t^N(x))^2\Big|\eta_0^N\right)\leq 2(\eta_0^N(x))^2
+2\left(t\frac{C_1}{N}\sum_{y=1}^N\eta_0^N(y)+\frac{t^2C_1^2}{N^2}\left(\sum_{y=1}^N\eta_0^N(y)\right)^2\right)
\]
and then, under Assumption (A),
\[
\mathbb{E}\left((\eta_t^N(x))^2\right)\leq 2(\|\psi\|_\infty^2+\|\psi\|_\infty)+2tC_1\|\psi\|_\infty
+2\frac{t^2C_1^2}{N^2}\left(N^2\|\psi\|_\infty^2+N\|\psi\|_\infty\right).
\]
Consequently, Equation \eqref{equ 5.1 two} holds.

Utilizing Markov inequality, for any $\theta>0$,
\begin{align*}
\mathbb{P}\left(\eta_t^N(x)=l\Big|\eta_0^N\right)&\leq e^{-\theta l}\mathbb{E}\left(e^{\theta\eta_t^N(x)}\Big|\eta_0^N\right) \\
&\leq e^{-\theta l}e^{\theta\eta_0^N(x)}\mathbb{E}\left(e^{\theta\mathcal{Y}\left(t\frac{C_1}{N}\sum_{y=1}^N\eta_0^N(y)\right)}\Big|\eta_0^N\right)\\
&=e^{-\theta l}e^{\theta\eta_0^N(x)}e^{(e^\theta-1)(t\frac{C_1}{N}\sum_{y=1}^N\eta_0^N(y))}.
\end{align*}
Therefore, under Assumption (A),
\begin{align*}
\mathbb{P}\left(\eta_t^N(x)=l\right)&\leq e^{-\theta l}\mathbb{E}e^{\eta_0^N(x)(C_8^N(\theta, t)+\theta)}\mathbb{E}e^{\sum_{y\neq x}\eta_0^N(y)C_8^N(\theta, t)}\\
&=e^{-\theta l}e^{\psi(x/N)\left(e^{C_8^N(\theta, t)+\theta}-1\right)}\prod_{y\neq x}e^{\psi(y/N)\left(e^{C_8^N(\theta,t)}-1\right)}\\
&\leq e^{-\theta l}\exp\left\{\|\psi\|_\infty\left(e^{C_8^N(\theta, t)+\theta}-1+N\left(e^{C_8^N(\theta,t)}-1\right)\right)\right\},
\end{align*}
where $C_8^N(\theta, t)=\frac{tC_1}{N}(e^\theta-1)$. Since
\[
\lim_{N\rightarrow\infty}C_8^{N}(\theta, t)=0 \text{~and~}\lim_{N\rightarrow+\infty}N\left(e^{C_8^N(\theta,t)}-1\right)=tC_1(e^\theta-1),
\]
there exists $C_{9}=C_9(T, \theta)$ independent of $N, l$ such that
\[
\mathbb{P}\left(\eta_t^N(x)=l\right)\leq e^{-\theta l}e^{\|\psi\|_\infty C_{9}(\theta, T)}
\]
for all $0\leq t\leq T$, $N\geq 1, l\geq 0$ and $1\leq x\leq N$. Since $\theta$ is arbitrary in the last inequality, Equation \eqref{equ 5.2} holds and the proof is complete.

\qed

In this section, as in Section \ref{section four proof of hydro finite}, we still define
\[
\varrho_{t,k}^N(u)=\mathbb{P}\left(\eta_t^N(i)=k\right)
\]
for any $k\geq 0$ and $\frac{i-1}{N}<u\leq \frac{i}{N}$ for some $i\in \{1,\ldots, N\}$, then $\{\varrho_{t,k}^N(u)\}_{u\in \mathbb{T}}\in \mathcal{D}(\mathbb{T})$ for any $t\geq 0$, where $\mathcal{D}(\mathbb{T})$ is the set of c\`{a}dl\`{a}g functions from $\mathbb{T}$ to $\mathbb{R}$ endowed with the Skorohod topology.

\begin{lemma}\label{lemma 5.2}
For any given $T>0$ and $k\geq 0$, $\left\{\varrho^N_{t,k}: 0\leq t\leq T\right\}_{N\geq 1}$ are relatively compact with respect to the uniform topology of $C([0, T], \mathcal{D}(\mathbb{T}))$, where $C([0, T], \mathcal{D}(\mathbb{T}))$ is the set of continuous functions from $[0, T]$ to $\mathcal{D}(\mathbb{T})$.
\end{lemma}

\proof

Using Arzela-Ascoli lemma, we only need to show that $\left\{\varrho^N_{t,k}: 0\leq t\leq T\right\}_{N\geq 1}$ are uniformly bounded and equicontinuous. The uniform boundness is obvious since a probability is at most $1$. Now we check the equicontinuity. Here we only discuss the case where $k\geq 1$, since the case where $k=0$ can be checked similarly. Using Kolmogorov-Chapman equation, we have
\begin{align}\label{equ 5.3}
\frac{d}{dt}\varrho_{t,k}^N(u)=&-\frac{1}{N}\sum_{l=0}^{+\infty}\sum_{j=1}^N\phi_{k,l}(i/N, j/N)\mathbb{P}\left(\eta_t^N(i)=k, \eta_t^N(j)=l\right)\notag\\
&-\frac{1}{N}\sum_{l=1}^{+\infty}\sum_{j=1}^N\phi_{l,k}(j/N,i/N)\mathbb{P}\left(\eta_t^N(i)=k, \eta_t^N(j)=l\right)\\
&+\frac{1}{N}\sum_{l=1}^{+\infty}\sum_{j=1}^N\phi_{l,k-1}(j/N,i/N)\mathbb{P}\left(\eta_t^N(i)=k-1, \eta_t^N(j)=l\right)\notag\\
&+\frac{1}{N}\sum_{l=0}^{+\infty}\sum_{j=1}^N\phi_{k+1,l}(i/N,j/N)\mathbb{P}\left(\eta_t^N(i)=k+1, \eta_t^N(j)=l\right).\notag
\end{align}
for any $(i-1)/N<u\leq i/N$. By Assumption (B),
\begin{align*}
\frac{1}{N}\sum_{l=0}^{+\infty}\sum_{j=1}^N\phi_{k,l}(i/N, j/N)\mathbb{P}\left(\eta_t^N(i)=k, \eta_t^N(j)=l\right)
& \leq \frac{1}{N}\sum_{l=0}^{+\infty}\sum_{j=1}^N\phi_{k,l}(i/N, j/N)\mathbb{P}\left(\eta_t^N(j)=l\right)\\
&\leq \frac{1}{N}\sum_{l=0}^{+\infty}\sum_{j=1}^NC_1k\mathbb{P}\left(\eta_t^N(j)=l\right)\leq C_1k.
\end{align*}
Similarly,
\[
\frac{1}{N}\sum_{l=0}^{+\infty}\sum_{j=1}^N\phi_{k+1,l}(i/N,j/N)\mathbb{P}\left(\eta_t^N(i)=k+1, \eta_t^N(j)=l\right)\leq (k+1)C_1.
\]
By Assumption (B) and Lemma \ref{lemma 6.1}, for $0\leq t\leq T$,
\begin{align*}
&\frac{1}{N}\sum_{l=1}^{+\infty}\sum_{j=1}^N\phi_{l,k}(j/N,i/N)\mathbb{P}\left(\eta_t^N(i)=k, \eta_t^N(j)=l\right)\\
& \leq \frac{1}{N}\sum_{l=1}^{+\infty}\sum_{j=1}^NC_1l\mathbb{P}\left(\eta_t^N(j)=l\right)\\
&=\frac{1}{N}\sum_{j=1}^NC_1\mathbb{E}\eta_t^N(j)\leq C_1(TC_1+1)\|\psi\|_\infty.
\end{align*}
Similarly, for $0\leq t\leq T$,
\[
\frac{1}{N}\sum_{l=1}^{+\infty}\sum_{j=1}^N\phi_{l,k-1}(j/N,i/N)\mathbb{P}\left(\eta_t^N(i)=k-1, \eta_t^N(j)=l\right)
\leq C_1(TC_1+1)\|\psi\|_\infty.
\]
Hence, by Equation \eqref{equ 5.3}, for $0\leq t\leq T$ and $u\in \mathbb{T}$,
\begin{equation}\label{equ 5.3 two}
|\frac{d}{dt}\varrho_{t,k}^N(u)|\leq (2k+1)C_1+2C_1(TC_1+1)\|\psi\|_\infty
\end{equation}
and therefore
\[
|\varrho_{t,k}^N(u)-\varrho_{s,k}^N(u)|\leq \left((2k+1)C_1+2C_1(TC_1+1)\|\psi\|_\infty\right)|t-s|
\]
for any $0\leq s,t\leq T$. As a result, $\left\{\varrho^N_{t,k}: 0\leq t\leq T\right\}_{N\geq 1}$ are equicontinuous and the proof is complete.

\qed

\begin{lemma}\label{lemma 5.3}
For any given $T>0$, the solution to Equation \eqref{equ 2.3 mean field ODE infinite types} under constraints \eqref{equ constraints one} and \eqref{equ constraints two} exists for $0\leq t\leq T$.
\end{lemma}

\proof

Using Lemma \ref{lemma 5.2}, we let $\{\varrho_{t,k}^\infty: 0\leq t\leq T\}$ be the limit in $\mathcal{D}(\mathbb{T})$ of a subsequence 
\[
\left\{\varrho^{N_j}_{t,k}: 0\leq t\leq T\right\}_{j\geq 1}
\]
of $\left\{\varrho^{N}_{t,k}: 0\leq t\leq T\right\}_{N\geq 1}$ for all $k\geq 0$. For simplicity, we still write this subsequence as $\left\{\varrho^N_{t,k}: 0\leq t\leq T\right\}_{N\geq 1}$. Now we check that $\{\{\varrho_{t,k}^\infty: 0\leq t\leq T\}\}_{k\geq 0}$ is a solution to Equation \eqref{equ 2.3 mean field ODE infinite types} under constraints \eqref{equ constraints one} and \eqref{equ constraints two} for $0\leq t\leq T$.

Using Lemma \ref{lemma 3.2 approximated independence infinite K},
\[
\frac{1}{N}\sum_{l=0}^{M}\sum_{j=1}^N\phi_{k,l}(i/N, j/N)|\mathbb{P}\left(\eta_t^N(i)=k, \eta_t^N(j)=l\right)-\mathbb{P}\left(\eta_t^N(i)=k\right) \mathbb{P}\left(\eta_t^N(j)=l\right)|=0
\]
for any $M\geq 1$. Then, using Equation \eqref{equ 5.2}, we have
\[
\frac{1}{N}\sum_{l=0}^{+\infty}\sum_{j=1}^N\phi_{k,l}(i/N, j/N)|\mathbb{P}\left(\eta_t^N(i)=k, \eta_t^N(j)=l\right)-\mathbb{P}\left(\eta_t^N(i)=k\right) \mathbb{P}\left(\eta_t^N(j)=l\right)|=0.
\]
Similar equations hold for the other three terms in the right-hand side of Equation \eqref{equ 5.3}. Therefore, let $N\rightarrow+\infty$ in Equation \eqref{equ 5.3}, we have
\begin{align}\label{equ 5.4}
\lim_{N\rightarrow+\infty}\frac{d}{dt}\varrho_{t,k}^\infty(u)=&-\int_\mathbb{T}\sum_{l=0}^{+\infty}\phi_{k,l}(u, v)\varrho_{t,k}^\infty(u)\varrho_{t,l}^\infty(v)dv\notag\\
&-\int_{\mathbb{T}}\sum_{l=1}^{+\infty}\phi_{l,k}(v,u)\varrho_{t,k}^\infty(u)\varrho_{t,l}^\infty(v)dv\\
&+\int_{\mathbb{T}}\sum_{l=1}^{+\infty}\phi_{l,k-1}(v,u)\varrho_{t,k-1}^\infty(u)\varrho_{t,l}^\infty(v)dv\notag\\
&+\int_{\mathbb{T}}\sum_{l=0}^{+\infty}\phi_{k+1,l}(u,v)\varrho_{t,k+1}^\infty(u)\varrho_{t,l}^\infty(v)dv \notag
\end{align}
for $k\geq 1$. Similarly,
\begin{equation}\label{equ 5.5}
\lim_{N\rightarrow+\infty}\frac{d}{dt}\varrho_{t,0}^N(u)=\sum_{l=0}^\infty\int_{\mathbb{T}}\phi_{1,l}(u,v)\varrho_{t,1}^{\infty}(u)\varrho_{t,l}^\infty(v)dv
-\sum_{l=1}^\infty\int_{\mathbb{T}}\phi_{l,0}(v,u)\varrho_{t,0}^\infty(u)\varrho_{t,l}^\infty(v)dv.
\end{equation}
According to Assumption (A), for all $k\geq 0$,
\begin{equation}\label{equ 5.6}
\varrho_{0,k}^\infty(u)=e^{-\psi(u)}\frac{\psi(u)^k}{k!}.
\end{equation}
By Equation \eqref{equ 5.3 two} and the dominated convergence theorem,
\begin{align*}
\varrho_{t,k}^\infty(u)&=\lim_{N\rightarrow+\infty}\varrho_{t,k}^N(u) =\lim_{N\rightarrow+\infty}\left(\varrho_{0,k}^N(u)+\int_0^t\frac{d}{ds}\varrho_{s,k}^N(u)ds\right)\\
&=\varrho_{0,k}^{\infty}(u)+\int_0^t\lim_{N\rightarrow+\infty}\frac{d}{ds}\varrho_{s,k}^N(u)ds.
\end{align*}
Therefore, by Equations \eqref{equ 5.4}-\eqref{equ 5.6}, $\{\{\varrho_{t,k}^\infty: 0\leq t\leq T\}\}_{k\geq 0}$ is a solution to Equation \eqref{equ 2.3 mean field ODE infinite types} for $0\leq t\leq T$. Lemma \ref{lemma 5.1} ensures that $\{\{\varrho_{t,k}^\infty: 0\leq t\leq T\}\}_{k\geq 0}$ satisfies constraints \eqref{equ constraints one}, \eqref{equ constraints two} and the proof is complete. 

\qed

\begin{lemma}\label{lemma 5.4}
For any given $T>0$, the solution to Equation \eqref{equ O-U} under constraints \eqref{equ constraints one} and \eqref{equ constraints two} for $0\leq t\leq T$ is unique.
\end{lemma}

\proof

Assuming that $\{\rho_{t,k}^\infty: 0\leq t\leq T\}_{k\geq 0}$ and $\{\widehat{\rho}_{t,k}^\infty: 0\leq t\leq T\}_{k\geq 0}$ are both solutions to Equation \eqref{equ O-U} under constraints \eqref{equ constraints one} and \eqref{equ constraints two}. Now only need to show that $\rho_{t,k}^\infty(u)=\widehat{\rho}_{t,k}^\infty(u)$ for any $k\geq 0$ and $u\in \mathbb{T}$.

For any $M\geq 1$, we define
\[
\gamma_{M,t}=\max\left\{\sup_{1\leq k\leq M, u\in \mathbb{T}}k\left|\rho_{t,k}^\infty(u)-\widehat{\rho}_{t,k}^\infty(u)\right|, \text{~}
\sup_{u\in \mathbb{T}}\left|\rho_{t,0}^\infty(u)-\widehat{\rho}_{t,0}^\infty(u)\right|\right\},
\]
\[
C_{10}(M)=\sup_{0\leq t\leq T, u\in \mathbb{T}}\left(\sum_{l=M+1}^{+\infty}l\rho_{t,l}^\infty(u)\right), \text{\quad}
\widehat{C}_{10}(M)=\sup_{0\leq t\leq T, u\in \mathbb{T}}\left(\sum_{l=M+1}^{+\infty}l\widehat{\rho}_{t,l}^\infty(u)\right),
\]
\[
C_{11}=\sup_{0\leq t\leq T, u\in \mathbb{T}}\sum_{l=0}^{+\infty}l\rho_{t,l}^\infty(u), \text{\quad}
\widehat{C}_{11}=\sup_{0\leq t\leq T, u\in \mathbb{T}}\sum_{l=0}^{+\infty}l\widehat{\rho}_{t,l}^\infty(u)
\]
and
\[
C_{12}=\sup_{0\leq t\leq T, u\in \mathbb{T}}\sum_{l=0}^{+\infty}l^2\rho_{t,l}^\infty(u)\text{\quad and \quad}
\widehat{C}_{12}=\sup_{0\leq t\leq T, u\in \mathbb{T}}\sum_{l=0}^{+\infty}l^2\widehat{\rho}_{t,l}^\infty(u).
\]
Using Assumption (B) and triangle inequality, we have
\begin{align*}
&k\left|\frac{d}{dt}\rho_{t,k}^\infty(u)-\frac{d}{dt}\widehat{\rho}_{t,k}^\infty(u)\right|\\
&\leq C_1M\gamma_{M,t}+\widehat{C}_{12}C_1M\gamma_{M,t}+C_1\widehat{C}_{12}\left(C_{10}(M)+\widehat{C}_{10}(M)\right)\\
&\text{\quad}+C_1\gamma_{M,t}C_{11}+C_1\widehat{C}_{11}M\gamma_{M,t}+C_1\widehat{C}_{11}\left(C_{10}(M)+\widehat{C}_{10}(M)\right)\\
&\text{\quad}+2C_1\gamma_{M,t}C_{11}+C_1(\widehat{C}_{11}+1)M\gamma_{M,t}+C_1(\widehat{C}_{11}+1)\left(C_{10}(M)+\widehat{C}_{10}(M)\right)\\
&\text{\quad}+C_1M\gamma_{M,t}+\widehat{C}_{12}C_1M\gamma_{M,t}+C_1\widehat{C}_{12}\left(C_{10}(M)+\widehat{C}_{10}(M)\right)\\
&\text{\quad}+C_1M\left(C_{10}(M)+\widehat{C}_{10}(M)\right)
\end{align*}
for any $k\leq M, u\in \mathbb{T}$ and $0\leq t\leq T$. Then, using Grownwall's inequality, we have
\[
\gamma_{M,t}\leq C_{13}(M)e^{C_{14}(M)t}
\]
for $0\leq t\leq T$, where
\[
C_{13}(M)=T\left(C_{10}(M)+\widehat{C}_{10}(M)\right)\left(2C_1\widehat{C}_{12}+C_1\widehat{C}_{11}+C_1(\widehat{C}_{11}+1)+C_1M\right)
\]
and
\[
C_{14}(M)=M\left(2C_1+2\widehat{C}_{12}C_1+C_1\widehat{C}_1+C_1(\widehat{C}_{11}+1)\right)
+3C_1C_{11}.
\]
Since $C_{14}(M)=O(M)$, by constraint \eqref{equ constraints two}, $\lim_{M\rightarrow+\infty}C_{13}(M)e^{C_{14}(M)t}=0$. Since $\gamma_{M,t}$ is increasing with $M$, we have $\gamma_{M,t}=0$ for all $M\geq 1$ and the proof is complete.

\qed

At last, we prove Theorem \ref{theorem 2.2 mean field limit infinite}.

\proof[Proof of Theorem \ref{theorem 2.2 mean field limit infinite}]

By Lemmas \ref{lemma 5.3} and \ref{lemma 5.4}, the solution $\{\rho_{t,k}^\infty\}_{k\geq 0}$ to Equation \eqref{equ O-U} under constraints \eqref{equ constraints one} and \eqref{equ constraints two} exists for $t\geq 0$ and is unique. Furthermore, according to the proof of Lemma \ref{lemma 5.3},
\[
\lim_{N\rightarrow+\infty}\sup_{u\in \mathbb{T}, 0\leq t\leq T}\left|\varrho_{t,k}^N(u)-\rho_{t,k}^\infty(u)\right|=0
\]
for any given $T>0$ and $k\geq 0$. As a result,
\[
\lim_{N\rightarrow+\infty}\mathbb{E}\mu_{t,k}^N(f)=\int_{\mathbb{T}}\rho_{t,k}^\infty(u)f(u)du
\]
for any $t\geq 0, k\geq 0$ and $f\in C(\mathbb{T})$. Hence, to prove Theorem \ref{theorem 2.2 mean field limit infinite}, we only need to show that
\[
\lim_{N\rightarrow+\infty}{\rm Var}\left(\mu_{t,k}^N(f)\right)=0.
\]
Using Lemma \ref{lemma 3.2 approximated independence infinite K} and an analysis similar with that leading to Equation \eqref{equ 4.1}, we have
\[
{\rm Var}\left(\mu_{t,k}^N(f)\right)=O(N^{-0.9})
\]
and the proof is complete.

\qed

\section{The proof of Theorem \ref{theorem 2.3 fluctuation limit}}\label{section six proof of fluctuation}
In this section, we prove Theorem \ref{theorem 2.3 fluctuation limit}. For later use, we first introduce some notations. For any $k\in \{0,\ldots, \mathcal{K}\}$, $t\geq 0$ and $x\in \{1,\ldots,N\}$, we define
\[
\widehat{1}_{\{\eta_t^N(x)=k\}}=1_{\{\eta_t^N(x)=k\}}-P(\eta_t^N(x)=k).
\]
For any $k_1, k_2\in \{0,\ldots,\mathcal{K}\}$, $t\geq 0$ and $h\in C(\mathbb{T}^2)$, we define
\[
\vartheta_{t,k_1, k_2}^N(h)=\frac{1}{N^{1.5}}\sum_{x=1}^N\sum_{y=1}^N\widehat{1}_{\{\eta_t^N(x)=k_1\}}\widehat{1}_{\{\eta_t^N(y)=k_2\}}h(x/N, y/N).
\]

As preliminaries, we give two lemmas.

\begin{lemma}\label{lemma 6.1}
For any given $T>0$, $k_1, k_2\in \{0,\ldots, \mathcal{K}\}$ and $h\in C(\mathbb{T}^2)$,
\[
\lim_{N\rightarrow+\infty}\sup_{0\leq t\leq T}\mathbb{E}\left|\vartheta^N_{t,k_1, k_2}(h)\right|=0.
\]
\end{lemma}

\begin{lemma}\label{lemma 6.2}
For any given $T>0$, $\left\{\{V_{t,k}^N:0\leq t\leq T\}_{k=0}^{\mathcal{K}}\right\}_{N\geq 1}$ are tight with respect to the topology of $D\left([0, T], \left(\left(C^\infty(\mathbb{T})\right)^\prime\right)^{\mathcal{K}+1}\right)$.
\end{lemma}

\proof[Proof of Lemma \ref{lemma 6.1}]

Let $\mathcal{H}$ be the set of $h\in C(\mathbb{T}^2)$ with the form $h(u,v)=f(u)g(v)$ for all $u,v\in \mathbb{T}$ and some $f,g\in C(\mathbb{T})$. For $h\in \mathcal{H}$,
\[
\vartheta_{t, k_1, k_2}^N(h)=\frac{1}{\sqrt{N}}V_{t,k_1}^N(f)V_{t,k_2}^N(g)
\]
for some $f,g\in C(\mathbb{T})$. According to the analysis leading to Equation \eqref{equ 4.1}, for any $0\leq k\leq \mathcal{K}$ and $f\in C(\mathbb{T})$, we have
\begin{align*}
\sup_{0\leq t\leq T}\mathbb{E}\left((V_{t,k}^N(f))^2\right)&=\sup_{0\leq t\leq T}{\rm Var}(V_{t,k}^N(f))\\
&=\sup_{0\leq t\leq T}{\rm Var}(\sqrt{N}\mu_{t,k}^N(f))=NO(N^{-1})=O(1).
\end{align*}
Therefore, using Cauchy-Schwartz inequality,
\begin{equation}\label{equ 6.1}
\sup_{0\leq t\leq T}\mathbb{E}\left|\vartheta^N_{t,k_1, k_2}(h)\right|=O(1)N^{-\frac{1}{2}}
\end{equation}
for all $h\in \mathcal{H}$ and hence for all $h\in {\rm span}(\mathcal{H})$. Note that $O(1)$ relies on the choice of $h\in {\rm span}(\mathcal{H})$ in Equation \eqref{equ 6.1}.

For general $h\in C(\mathbb{T}^2)$,
\[
\vartheta_{t,k_1,k_2}^N(h)=\frac{1}{N^{1.5}}\sum_{x=1}^N\widehat{1}_{\{\eta_t^N(x)=k_1\}}Z_{x,k_2,t}^N(h),
\]
where
\[
Z_{x,k_2,t}^N(h)=\sum_{y=1}^N\widehat{1}_{\{\eta_t^N(y)=k_2\}}h(x/N,y/N).
\]
Hence, using Cauchy-Schwartz inequality twice,
\begin{align*}
\mathbb{E}|\vartheta_{t,k_1,k_2}^N(h)|&\leq \frac{1}{N^{1.5}}\mathbb{E}\left(\sqrt{\sum_{x=1}^N\widehat{1}^2_{\{\eta_t^N(x)=k_1\}}}\sqrt{\sum_{x=1}^N(Z_{x,k_2,t}^N(h))^2}\right)\\
&\leq \frac{1}{N^{1.5}}\sqrt{\mathbb{E}\sum_{x=1}^N\widehat{1}^2_{\{\eta_t^N(x)=k_1\}}}\sqrt{\mathbb{E}\sum_{x=1}^N(Z_{x,k_2,t}^N(h))^2}\\
&\leq \frac{1}{N^{1.5}}\sqrt{4N}\sqrt{\mathbb{E}\sum_{x=1}^N(Z_{x,k_2,t}^N(h))^2}.
\end{align*}
According to the analysis leading to Equation \eqref{equ 4.1},
\[
\mathbb{E}((Z_{x,k_2,t}^N(h))^2)\leq (2+C_2)N\|h\|_\infty^2,
\]
where $C_2=C_2(T)$ defined as in Lemma \ref{lemma 3.1 approximated independence finite K}. As a result, for general $h\in C(\mathbb{T}^2)$,
\begin{align}\label{equ 6.2}
\mathbb{E}|\vartheta_{t,k_1,k_2}^N(h)|\leq \frac{1}{N^{1.5}}\sqrt{4N}\sqrt{N^2(2+C_2)\|h\|_\infty^2}=2\sqrt{2+C_2}\|h\|_\infty.
\end{align}
For given $h\in C(\mathbb{T}^2)$ and any $\epsilon>0$, since ${\rm span}(\mathcal{H})$ is dense in $C(\mathbb{T}^2)$, there exists $h_\epsilon\in {\rm span}(\mathcal{H})$ such that $\|h-h_\epsilon\|\leq \epsilon$. Hence, by Equation \eqref{equ 6.2},
\[
\mathbb{E}|\vartheta_{t,k_1,k_2}^N(h)|\leq \mathbb{E}|\vartheta_{t,k_1,k_2}^N(h_\epsilon)|+2\sqrt{2+C_2}\|h-h_\epsilon\|_\infty
\leq \mathbb{E}|\vartheta_{t,k_1,k_2}^N(h_\epsilon)|+2\sqrt{2+C_2}\epsilon.
\]
Since $h_\epsilon\in {\rm span}(\mathcal{H})$, let $N\rightarrow+\infty$ in the last inequality, we have
\[
\limsup_{N\rightarrow+\infty}\sup_{0\leq t\leq T}\mathbb{E}\left|\vartheta^N_{t,k_1, k_2}(h)\right|\leq 2\sqrt{2+C_2}\epsilon
\]
by utilizing Equation \eqref{equ 6.1}. Since $\epsilon$ is arbitrary, let $\epsilon\rightarrow 0$ in the last inequality and then the proof is complete.

\qed

\proof[Proof of Lemma \ref{lemma 6.2}]

Utilizing Aldous' criteria, we only need to show that
\begin{equation}\label{equ 6.3}
\lim_{M\rightarrow+\infty}\limsup_{N\rightarrow+\infty}\mathbb{P}\left(|V_{t,k}^N(f)|\geq M\right)=0
\end{equation}
for any $t\geq 0, f\in C(\mathbb{T}), 0\leq k\leq \mathcal{K}$ and
\begin{equation}\label{equ 6.4}
\lim_{\delta\rightarrow 0}\limsup_{N\rightarrow+\infty}\sup_{\sigma\in \mathcal{T}, s\leq\delta}\mathbb{P}\big(|V_{\sigma+s,k}^N(f)-V_{\sigma,k}^N(f)|>\epsilon\big)=0
\end{equation}
for any $\epsilon>0, f\in C(\mathbb{T}), 0\leq k\leq \mathcal{K}$, where $\mathcal{T}$ is the set of stopping times of $\{\eta_t^N\}_{t\geq 0}$ bounded by $T$.

As we have shown in the proof of Lemma \ref{lemma 6.1}, $\mathbb{E}\left((V_{t,k}^N(f))^2\right)=O(1)$ and hence Equation \eqref{equ 6.3} follows from Markov inequality. We only need to check Equation \eqref{equ 6.4}.

By Dynkin's martingale formula, let
\[
\mathcal{M}_{t,k}^N(f)=V_{t,k}^N(f)-\int_0^t(\mathcal{L}_N+\partial_s)V_{s,k}^N(f)ds,
\]
then $\{\mathcal{M}_{t,k}^N(f)\}_{t\geq 0}$ is a martingale with quadratic variation process $\{\langle\mathcal{M}_k^N(f)\rangle_t\}_{t\geq 0}$ given by
\[
\langle\mathcal{M}_k^N(f)\rangle_t=\int_0^t\Bigg(\mathcal{L}_N\left((V_{s,k}^N(f))^2\right)-2V_{s,k}^N(f)\mathcal{L}_NV_{s,k}^N(f)\Bigg)ds.
\]
To check Equation \eqref{equ 6.4}, we only need to show that
\begin{equation}\label{equ 6.5}
\lim_{\delta\rightarrow 0}\limsup_{N\rightarrow+\infty}\sup_{\sigma\in \mathcal{T}, s\leq\delta}\mathbb{P}\big(|\mathcal{M}_{\sigma+s,k}^N(f)-\mathcal{M}_{\sigma,k}^N(f)|>\epsilon\big)=0
\end{equation}
and
\begin{equation}\label{equ 6.6}
\lim_{\delta\rightarrow 0}\limsup_{N\rightarrow+\infty}\sup_{\sigma\in \mathcal{T}, s\leq\delta}\mathbb{P}\left(\left|\int_\sigma^{\sigma+s}(\mathcal{L}_N+\partial_u)V_{u,k}^N(f)du\right|>\epsilon\right)=0
\end{equation}
for any $\epsilon>0, f\in C(\mathbb{T}), 0\leq k\leq \mathcal{K}$.

We first check Equation \eqref{equ 6.5}. By Doob's inequality,
\begin{equation}\label{equ 6.7}
\mathbb{P}\big(|\mathcal{M}_{\sigma+s,k}^N(f)-\mathcal{M}_{\sigma,k}^N(f)|>\epsilon\big)\leq \frac{1}{\epsilon^2}\mathbb{E}\left(\langle\mathcal{M}_k^N(f)\rangle_{\sigma+s}-\langle\mathcal{M}_k^N(f)\rangle_{\sigma}\right).
\end{equation}
According to direct calculation,
\begin{align*}
&\mathcal{L}_N\left((V_{s,k}^N(f))^2\right)-2V_{s,k}^N(f)\mathcal{L}_NV_{s,k}^N(f)\\
&=\frac{1}{N^2}\sum_{x=1}^N\sum_{y\neq x}\sum_{m=0}^\mathcal{K}\sum_{l=0}^\mathcal{K}\phi_{m,l}(x/N, y/N)1_{\{\eta_s^N(x)=m\}}1_{\{\eta_s^N(y)=l\}}\left(q_{m,l,x,y}^k(f)\right)^2,
\end{align*}
where
\begin{align*}
q_{m,l,x,y}^k(f)=
\begin{cases}
f(y/N)-f(x/N) & \text{~if~}m=k \text{~and~}l=k-1,\\
-f(y/N)-f(x/N) & \text{~if~}m=k \text{~and~}l=k,\\
-f(x/N) & \text{~if~}m=k \text{~and~}l\neq k, k-1,\\
f(x/N)-f(y/N) & \text{~if~}m=k+1 \text{~and~}l=k,\\
-f(y/N) & \text{~if~}m\neq k+1, k \text{~and~}l=k,\\
f(x/N)+f(y/N) & \text{~if~}m=k+1 \text{~and~}l=k-1,\\
f(y/N) & \text{~if~}m\neq k+1,k \text{~and~}l=k-1,\\
f(x/N) & \text{~if~}m=k+1 \text{~and~}l\neq k, k-1,\\
0 &\text{~else}.
\end{cases}
\end{align*}
As a result, since $s\leq \delta$,
\[
\mathbb{E}\left(\langle\mathcal{M}_k^N(f)\rangle_{\sigma+s}-\langle\mathcal{M}_k^N(f)\rangle_{\sigma}\right)
\leq 4\delta \mathcal{K}^2\sup_{0\leq m,l\leq \mathcal{K}}\|\phi_{m,l}\|_\infty\|f\|_\infty^2
\]
and hence Equation \eqref{equ 6.5} follows from Equation \eqref{equ 6.7}.

At last, we check Equation \eqref{equ 6.6}. By direct calculation,
\[
(\mathcal{L}_N+\partial_u)V_{u,k}^N(f)
=\frac{1}{N^{1.5}}\sum_{x=1}^N\sum_{y\neq x}\sum_{m=0}^\mathcal{K}\sum_{l=0}^\mathcal{K}\varsigma^N_{u,m,l}(x,y)q_{m,l,x,y}^k(f),
\]
where
\[
\varsigma^N_{u,m,l}(x,y)=\phi_{m,l}(x/N, y/N)\left(1_{\{\eta_u^N(x)=m, \eta_u^N(y)=l\}}-\mathbb{P}(\eta_u^N(x)=m, \eta_u^N(y)=l)\right).
\]
Hence, $(\mathcal{L}_N+\partial_u)V_{u,k}^N(f)={\rm \uppercase\expandafter{\romannumeral1}}_u+{\rm \uppercase\expandafter{\romannumeral2}}_u+{\rm \uppercase\expandafter{\romannumeral3}}_u+{\rm \uppercase\expandafter{\romannumeral4}}_u$, where
\[
{\rm \uppercase\expandafter{\romannumeral1}}_u
=\frac{1}{N^{1.5}}\sum_{x=1}^N\sum_{y\neq x}\sum_{m=0}^\mathcal{K}\sum_{l=0}^\mathcal{K}\widehat{1}_{\{\eta_u^N(x)=m\}}\mathbb{P}(\eta_u^N(y)=l)\phi_{m,l}(x/N, y/N)q_{m,l,x,y}^k(f),
\]
\[
{\rm \uppercase\expandafter{\romannumeral2}}_u
=\frac{1}{N^{1.5}}\sum_{x=1}^N\sum_{y\neq x}\sum_{m=0}^\mathcal{K}\sum_{l=0}^\mathcal{K}\widehat{1}_{\{\eta_u^N(y)=l\}}\mathbb{P}(\eta_u^N(x)=m)\phi_{m,l}(x/N, y/N)q_{m,l,x,y}^k(f),
\]
\[
{\rm \uppercase\expandafter{\romannumeral3}}_u
=\frac{1}{N^{1.5}}\sum_{x=1}^N\sum_{y\neq x}\sum_{m=0}^\mathcal{K}\sum_{l=0}^\mathcal{K}\widehat{1}_{\{\eta_u^N(x)=m\}}\widehat{1}_{\{\eta_u^N(y)=l\}}\phi_{m,l}(x/N, y/N)q_{m,l,x,y}^k(f)
\]
and
\[
{\rm \uppercase\expandafter{\romannumeral4}}_u
=\frac{-1}{N^{1.5}}\sum_{x=1}^N\sum_{y\neq x}\sum_{m=0}^\mathcal{K}\sum_{l=0}^\mathcal{K}{\rm Cov}\left(1_{\{\eta_u^N(x)=m\}}, 1_{\{\eta_u^N(y)=l\}}\right)\phi_{m,l}(x/N, y/N)q_{m,l,x,y}^k(f).
\]
By Lemma \ref{lemma 3.1 approximated independence finite K},
\[
\sup_{\sigma\leq u\leq \sigma+s}\left|{\rm \uppercase\expandafter{\romannumeral4}}_u\right|=\frac{N(N-1)}{N^{1.5}}O(N^{-1})=O(N^{-0.5}).
\]
Hence, to prove Equation \eqref{equ 6.6}, we only need to show that
\begin{equation}\label{equ 6.8}
\lim_{\delta\rightarrow 0}\limsup_{N\rightarrow+\infty}\sup_{\sigma\in \mathcal{T}, s\leq\delta}\mathbb{P}\left(\left|\int_\sigma^{\sigma+s}{\rm \uppercase\expandafter{\romannumeral1}}_udu\right|>\epsilon\right)=0,
\end{equation}
\begin{equation}\label{equ 6.9}
\lim_{\delta\rightarrow 0}\limsup_{N\rightarrow+\infty}\sup_{\sigma\in \mathcal{T}, s\leq\delta}\mathbb{P}\left(\left|\int_\sigma^{\sigma+s}{\rm \uppercase\expandafter{\romannumeral2}}_udu\right|>\epsilon\right)=0,
\end{equation}
and
\begin{equation}\label{equ 6.10}
\lim_{\delta\rightarrow 0}\limsup_{N\rightarrow+\infty}\sup_{\sigma\in \mathcal{T}, s\leq\delta}\mathbb{P}\left(\left|\int_\sigma^{\sigma+s}{\rm \uppercase\expandafter{\romannumeral3}}_udu\right|>\epsilon\right)=0.
\end{equation}
We first check Equation \eqref{equ 6.10}. According to the definition of $q_{m,l,x,y}^k(f)$, we have
\[
\phi_{m,l}(x/N, y/N)q_{m,l,x,y}^k(f)=h_{m,l,f}^k(x/N, y/N)
\]
for some $h_{m,l,f}^k\in C(\mathbb{T}^2)$. Hence,
\[
\int_\sigma^{\sigma+s}{\rm \uppercase\expandafter{\romannumeral3}}_udu=O(N^{-0.5})+\int_\sigma^{\sigma+s}\sum_{k=0}^\mathcal{K}\sum_{k=0}^\mathcal{K}\vartheta_{u,m,l}^N(h_{m,l,f}^k)du.
\]
Therefore, to prove Equation \eqref{equ 6.10}, we only need to show that
\begin{equation}\label{equ 6.11}
\lim_{\delta\rightarrow 0}\limsup_{N\rightarrow+\infty}\sup_{\sigma\in \mathcal{T}, s\leq\delta}\mathbb{P}\left(\left|\int_\sigma^{\sigma+s}\vartheta_{u,m,l}^N(h_{m,l,f}^k)du\right|>\epsilon\right)=0
\end{equation}
for any $0\leq m,l\leq \mathcal{K}$. Since
\[
\left|\int_\sigma^{\sigma+s}\vartheta_{u,m,l}^N(h_{m,l,f}^k)du\right|\leq \int_0^{T+\delta}|\vartheta_{u,m,l}^N(h_{m,l,f}^k)|du,
\]
Equation \eqref{equ 6.11} follows from Lemma \ref{lemma 6.1} and Markov inequality. Therefore, Equation \eqref{equ 6.10} holds.

 Finally, we only need to prove Equations \eqref{equ 6.8} and \eqref{equ 6.9}. Here we only check Equation \eqref{equ 6.8} since Equation \eqref{equ 6.9} can be checked in the same way. Using Lemma \ref{lemma 3.1 approximated independence finite K} and an analysis simliar with that leading to Equation \eqref{equ 4.5}, we have
\[
\sup_{0\leq u\leq T+1}\mathbb{E}\left({\rm \uppercase\expandafter{\romannumeral1}}_u^2\right)=O(1).
\]
Hence, by Cauchy-Schwartz inequality,
\begin{align}\label{equ 6.12}
\mathbb{E}\left(\left(\int_\sigma^{\sigma+s}{\rm \uppercase\expandafter{\romannumeral1}}_udu\right)^2\right)
&\leq \mathbb{E}\left(s\int_\sigma^{\sigma+s}{\rm \uppercase\expandafter{\romannumeral1}}_u^2du\right) \notag\\
&\leq \delta\int_0^{T+1}\mathbb{E}\left({\rm \uppercase\expandafter{\romannumeral1}}_u^2\right)du=\delta O(1)
\end{align}
for sufficiently small $\delta$. Consequently, Equation \eqref{equ 6.8} follows from Equation \eqref{equ 6.12} and Markov inequality. Hence, Equation \eqref{equ 6.6} holds and the proof is complete.

\qed

At last, we prove Theorem \ref{theorem 2.3 fluctuation limit}.

\proof[Proof of Theorem \ref{theorem 2.3 fluctuation limit}]

By Lemma \ref{lemma 6.2}, we only need to show that any weak limit of a subsequence of $\{\{V_{t,k}^N: 0\leq t\leq T\}_{0\leq k\leq \mathcal{K}}\}_{N\geq 1}$ is a solution to the martingale problem with respect to Equation \eqref{equ O-U}. From now on we assume that $\{\widehat{V}_{t,k}: 0\leq t\leq T\}_{0\leq k\leq \mathcal{K}}$ is a weak limit of a subsequence of $\{\{V_{t,k}^N: 0\leq t\leq T\}_{0\leq k\leq \mathcal{K}}\}_{N\geq 1}$. For simplicity, we still write this subsequence as $\{\{V_{t,k}^N: 0\leq t\leq T\}_{0\leq k\leq \mathcal{K}}\}_{N\geq 1}$.

For any $G\in C_c^\infty(\mathbb{R}^{\mathcal{K}+1})$ and $f_0, f_1,\ldots,f_k\in C^\infty(\mathbb{T})$, by Dynkin's martingale formula,
\[
\left\{G\left(\{V_{t,k}^N(f_k)\}_{k=0}^\mathcal{K}\right)-G(\{V_{0,k}^N(f_k)\}_{k=0}^\mathcal{K})-\int_0^t(\partial_s+\mathcal{L}_N)
G(\{V_{s,k}^N(f_k)\}_{k=0}^\mathcal{K})ds\right\}_{0\leq t\leq T}
\]
is a martingale. By direct calculation,
\begin{align}\label{equ 6.13}
&\partial_sG(\{V_{s,k}^N(f_k)\}_{k=0}^\mathcal{K})=-\frac{1}{N^{1.5}}\sum_{k=0}^{\mathcal{K}}\sum_{m=0}^{\mathcal{K}}\sum_{l=0}^{\mathcal{K}}\sum_{x=1}^N\sum_{y\neq x}\Bigg(\partial_kG(\{V_{s,k}^N(f_k)\}_{k=0}^\mathcal{K}) \notag\\
&\text{\quad\quad}\times\mathbb{P}(\eta_s^N(x)=m, \eta_s^N(y)=l)\phi_{m,l}(x/N, y/N)q_{m,l,x,y}^k(f_k)\Bigg)
\end{align}
and
\begin{align}\label{equ 6.14}
\mathcal{L}_NG(\{V_{s,k}^N(f_k)\}_{k=0}^\mathcal{K})
&=\frac{1}{N}\sum_{x=1}^N\sum_{y\neq x}\sum_{m=0}^\mathcal{K}\sum_{l=0}^\mathcal{K}\Bigg(\phi_{m,l}(x/N, y/N)1_{\{\eta_s^N(x)=m, \eta_s^N(y)=l\}} \notag\\
&\times\Big(G\big(\{V_{s,k}^N(f_k)+\frac{q_{m,l,x,y}^k(f_k)}{\sqrt{N}}\}_{k=0}^\mathcal{K}\big)-G(\{V_{s,k}^N(f_k)\}_{k=0}^\mathcal{K})\Big)\Bigg).
\end{align}
In Equation \eqref{equ 6.14}, utilizing Taylor's expansion formula up to the second order with the Lagrange remainder, we have
\[
\mathcal{L}_NG(\{V_{s,k}^N(f_k)\}_{k=0}^\mathcal{K})={\rm \uppercase\expandafter{\romannumeral5}}^N_s+{\rm \uppercase\expandafter{\romannumeral6}}^N_s+\varepsilon_{7,s}^N,
\]
where
\begin{align*}
{\rm \uppercase\expandafter{\romannumeral5}}^N_s
&=\frac{1}{N^{1.5}}\sum_{k=0}^\mathcal{K}\sum_{x=1}^N\sum_{y\neq x}\sum_{m=0}^{\mathcal{K}}\sum_{l=0}^{\mathcal{K}}\Bigg(\partial_kG(\{V_{s,k}^N(f_k)\}_{k=0}^\mathcal{K})\phi_{m,l}(x/N, y/N)\\
&\text{\quad\quad}\times 1_{\{\eta_s^N(x)=m, \eta_s^N(y)=l\}}q_{m,l,x,y}^k(f_k)\Bigg),
\end{align*}
\begin{align*}
{\rm \uppercase\expandafter{\romannumeral6}}^N_s&=\frac{1}{2N^2}\sum_{k_1=0}^{\mathcal{K}}\sum_{k_2=0}^{\mathcal{K}}\sum_{x=1}^N\sum_{y\neq x}\sum_{m=0}^\mathcal{K}\sum_{l=0}^{\mathcal{K}}\Bigg(\phi_{m,l}(x/N, y/N)\partial^2_{k_1k_2}G(\{V_{s,k}^N(f_k)\}_{k=0}^\mathcal{K})\\
&\text{\quad\quad}\times 1_{\{\eta_s^N(x)=m, \eta_s^N(y)=l\}}q_{m,l,x,y}^{k_1}(f_{k_1})q_{m,l,x,y}^{k_2}(f_{k_2})\Bigg)
\end{align*}
and
\begin{equation}\label{equ 6.15}
\sup_{0\leq s\leq T}|\varepsilon_{7,s}^N|=O(N^{-0.5}).
\end{equation}
Furthermore, ${\rm \uppercase\expandafter{\romannumeral6}}^N_s
={\rm \uppercase\expandafter{\romannumeral8}}^N_s+\varepsilon_{9,s}^N$, where
\begin{align*}
{\rm \uppercase\expandafter{\romannumeral8}}^N_s&=\frac{1}{2N^2}\sum_{k_1=0}^{\mathcal{K}}\sum_{k_2=0}^{\mathcal{K}}\sum_{x=1}^N\sum_{y\neq x}\sum_{m=0}^\mathcal{K}\sum_{l=0}^{\mathcal{K}}\Bigg(\phi_{m,l}(x/N, y/N)\partial^2_{k_1k_2}G(\{V_{s,k}^N(f_k)\}_{k=0}^\mathcal{K})\\
&\text{\quad\quad}\times \mathbb{P}\left(\eta_s^N(x)=m, \eta_s^N(y)=l\right)q_{m,l,x,y}^{k_1}(f_{k_1})q_{m,l,x,y}^{k_2}(f_{k_2})\Bigg)
\end{align*}
and
\begin{align*}
\varepsilon_{9,s}^N&=\frac{1}{2N^2}\sum_{k_1=0}^{\mathcal{K}}\sum_{k_2=0}^{\mathcal{K}}\sum_{x=1}^N\sum_{y\neq x}\sum_{m=0}^\mathcal{K}\sum_{l=0}^{\mathcal{K}}\Bigg(\partial^2_{k_1k_2}G(\{V_{s,k}^N(f_k)\}_{k=0}^\mathcal{K})\\
&\text{\quad\quad}\times \varsigma^N_{s,m,l}(x,y)q_{m,l,x,y}^{k_1}(f_{k_1})q_{m,l,x,y}^{k_2}(f_{k_2})\Bigg).
\end{align*}
According to an analysis similar with that leading to Equation \eqref{equ 6.6}, we have
\[
\sup_{0\leq s\leq T}\mathbb{E}|N^{0.5}\varepsilon_{9,s}^N|=O(1)
\]
and hence
\begin{equation}\label{equ 6.16}
\sup_{0\leq s\leq T}\mathbb{E}|\varepsilon_{9,s}^N|=O(N^{-0.5}).
\end{equation}
Using Lemma \ref{lemma 3.1 approximated independence finite K}, we have  $\sup_{0\leq s\leq T}\left|{\rm \uppercase\expandafter{\romannumeral8}}^N_s-{\rm \uppercase\expandafter{\romannumeral10}}^N_s\right|=O(N^{-1})$, where
\begin{align*}
{\rm \uppercase\expandafter{\romannumeral10}}^N_s&=\frac{1}{2N^2}\sum_{k_1=0}^{\mathcal{K}}\sum_{k_2=0}^{\mathcal{K}}\sum_{x=1}^N\sum_{y\neq x}\sum_{m=0}^\mathcal{K}\sum_{l=0}^{\mathcal{K}}\Bigg(\phi_{m,l}(x/N, y/N)\partial^2_{k_1k_2}G(\{V_{s,k}^N(f_k)\}_{k=0}^\mathcal{K})\\
&\text{\quad\quad}\times \mathbb{P}\left(\eta_s^N(x)=m\right)\mathbb{P}\left(\eta_s^N(y)=l\right)q_{m,l,x,y}^{k_1}(f_{k_1})q_{m,l,x,y}^{k_2}(f_{k_2})\Bigg).
\end{align*}
By Lemma \ref{lemma 4.2}, when we replace $\mathbb{P}\left(\eta_s^N(x)=m\right)\mathbb{P}\left(\eta_s^N(y)=l\right)$ by $\rho^{\mathcal{K}}_{s,m}(x/N)\rho^{\mathcal{K}}_{s,l}(y/N)$ in ${\rm \uppercase\expandafter{\romannumeral10}}^N_s$, the error is $o(1)$. Furthermore, using the definitions of $q_{m,l,x,y}^{k}(f_{k}), b_{t,m,l}^k$, we have
\begin{align*}
&\phi_{m,l}(x/N, y/N)\rho^{\mathcal{K}}_{s,m}(x/N)\rho^{\mathcal{K}}_{s,l}(y/N)q_{m,l,x,y}^{k_1}(f_{k_1})q_{m,l,x,y}^{k_2}(f_{k_2})\\
&=(b_{s,m,l}^{k_1}f_{k_1})(x/N,y/N)(b_{s,m,l}^{k_2}f_{k_2})(x/N,y/N).
\end{align*}
Therefore,
\begin{align*}
{\rm \uppercase\expandafter{\romannumeral10}}^N_s=&\frac{1}{2}\sum_{k_1=0}^\mathcal{K}\sum_{k_2=0}^\mathcal{K}\sum_{m=0}^\mathcal{K}\sum_{l=0}^\mathcal{K}\int_{\mathbb{T}^2}\partial^2_{k_1k_2}
G(\{V_{s,k}^N(f_k)\}_{k=0}^\mathcal{K})(b_{s,m,l}^{k_1}f_{k_1}(u,v))(b_{s,m,l}^{k_2}f_{k_2}(u,v))dudv\\
&+\varepsilon^N_{11,s},
\end{align*}
where $\lim_{N\rightarrow+\infty}\sup_{0\leq s\leq T}|\varepsilon^N_{11,s}|=0$. Then, by Equation \eqref{equ 6.16}, we have
\begin{align}\label{equ 6.17}
{\rm \uppercase\expandafter{\romannumeral6}}^N_s=&\frac{1}{2}\sum_{k_1=0}^\mathcal{K}\sum_{k_2=0}^\mathcal{K}\sum_{m=0}^\mathcal{K}\sum_{l=0}^\mathcal{K}\int_{\mathbb{T}^2}\partial^2_{k_1k_2}
G(\{V_{s,k}^N(f_k)\}_{k=0}^\mathcal{K})(b_{s,m,l}^{k_1}f_{k_1}(u,v))(b_{s,m,l}^{k_2}f_{k_2}(u,v))dudv \notag\\
&+\varepsilon^N_{12,s},
\end{align}
where $\lim_{N\rightarrow+\infty}\sup_{0\leq s\leq T}\mathbb{E}|\varepsilon^N_{12,s}|=0$.

Utilizing Equation \eqref{equ 6.13} and the expression of $(\mathcal{L}_N+\partial_s)V_{s,k}^N(f)$ given in the proof of Lemma \ref{lemma 6.2}, we have
\begin{align*}
&\partial_sG(\{V_{s,k}^N(f_k)\}_{k=0}^\mathcal{K})+{\rm \uppercase\expandafter{\romannumeral5}}^N_s\\
&=\frac{1}{N^{1.5}}\sum_{k=0}^\mathcal{K}\sum_{x=1}^N\sum_{y\neq x}\sum_{m=0}^{\mathcal{K}}\sum_{l=0}^{\mathcal{K}}\Bigg(\partial_kG(\{V_{s,k}^N(f_k)\}_{k=0}^\mathcal{K})
\varsigma^N_{s,m,l}(x,y)q_{m,l,x,y}^k(f_k)\Bigg)\\
&=\sum_{k=0}^{\mathcal{K}}\partial_kG(\{V_{s,k}^N(f_k)\}_{k=0}^\mathcal{K})(\mathcal{L}_N+\partial_s)V_{s,k}^N(f_k).
\end{align*}
Let ${\rm \uppercase\expandafter{\romannumeral1}}_s, {\rm \uppercase\expandafter{\romannumeral2}}_s, {\rm \uppercase\expandafter{\romannumeral3}}_s, {\rm \uppercase\expandafter{\romannumeral4}}_s$ be defined as in the proof of Lemma \ref{lemma 6.2} and be rewritten as
\[
{\rm \uppercase\expandafter{\romannumeral1}}_s(f), {\rm \uppercase\expandafter{\romannumeral2}}_s(f), {\rm \uppercase\expandafter{\romannumeral3}}_s(f), {\rm \uppercase\expandafter{\romannumeral4}}_s(f)
 \]
to emphasize the dependence on $f$. According to the analysis given in the proof of Lemma \ref{lemma 6.2},
 \begin{align*}
 &\partial_sG(\{V_{s,k}^N(f_k)\}_{k=0}^\mathcal{K})+{\rm \uppercase\expandafter{\romannumeral5}}^N_s\\
 &=\sum_{k=0}^{\mathcal{K}}\partial_kG(\{V_{s,k}^N(f_k)\}_{k=0}^\mathcal{K})
 \left({\rm \uppercase\expandafter{\romannumeral1}}_s(f_k)+{\rm \uppercase\expandafter{\romannumeral2}}_s(f_k)+{\rm \uppercase\expandafter{\romannumeral3}}_s(f_k)+{\rm \uppercase\expandafter{\romannumeral4}}_s(f_k)\right)\\
 &=\sum_{k=0}^{\mathcal{K}}\partial_kG(\{V_{s,k}^N(f_k)\}_{k=0}^\mathcal{K})
 \left({\rm \uppercase\expandafter{\romannumeral1}}_s(f_k)+{\rm \uppercase\expandafter{\romannumeral2}}_s(f_k)\right)+\varepsilon_{13,s}^N,
 \end{align*}
where $\lim_{N\rightarrow+\infty}\sup_{0\leq s\leq T}\mathbb{E}|\varepsilon_{13,s}^N|=0$. By Lemma \ref{lemma 4.2},
\begin{equation}\label{equ 6.18}
\sup_{0\leq s\leq T}\left|\frac{1}{N}\sum_{y=1}^N\left(\mathbb{P}(\eta_s^N(y)=l)-\rho_{s,l}^\mathcal{K}(y/N)\right)\phi_{m,l}(x/N, y/N)q_{m,l,x,y}^k(f_k)\right|=o(1).
\end{equation}
By the definition of $\mathcal{P}_{s,m,l}^{k,1}f_k$,
\begin{equation}\label{equ 6.19}
\sup_{0\leq s\leq T, 1\leq x\leq N}\left|\frac{1}{N}\sum_{y=1}^N\rho_{s,l}^\mathcal{K}(y/N)\phi_{m,l}(x/N, y/N)q_{m,l,x,y}^k(f_k)-\mathcal{P}_{s,m,l}^{k,1}f_k(x/N)\right|=o(1).
\end{equation}
Using Equations \eqref{equ 6.18}, \eqref{equ 6.19}, Lemma \ref{lemma 3.1 approximated independence finite K} and Markov inequality, we have
\[
{\rm \uppercase\expandafter{\romannumeral1}}_s(f_k)=\sum_{m=0}^{\mathcal{K}}\sum_{l=0}^{\mathcal{K}}V_{s,m}^N(\mathcal{P}_{s,m,l}^{k,1}f_k)+\varepsilon_{14,s}^N,
\]
where $\lim_{N\rightarrow+\infty}\sup_{0\leq s\leq T}\mathbb{E}|\varepsilon_{14,s}^N|=0$. Similarly,
\[
{\rm \uppercase\expandafter{\romannumeral2}}_s(f_k)=\sum_{m=0}^{\mathcal{K}}\sum_{l=0}^{\mathcal{K}}V_{s,l}^N(\mathcal{P}_{s,m,l}^{k,2}f_k)+\varepsilon_{15,s}^N,
\]
where $\lim_{N\rightarrow+\infty}\sup_{0\leq s\leq T}\mathbb{E}|\varepsilon_{15,s}^N|=0$. Therefore,
\begin{align}\label{equ 6.20}
&\partial_sG(\{V_{s,k}^N(f_k)\}_{k=0}^\mathcal{K})+{\rm \uppercase\expandafter{\romannumeral5}}^N_s\\
&=\sum_{k=0}^{\mathcal{K}}\sum_{m=0}^{\mathcal{K}}\sum_{l=0}^{\mathcal{K}}\partial_kG(\{V_{s,k}^N(f_k)\}_{k=0}^\mathcal{K})
\left(V_{s,m}^N(\mathcal{P}_{s,m,l}^{k,1}f_k)+V_{s,l}^N(\mathcal{P}_{s,m,l}^{k,2}f_k)\right)+\varepsilon_{16,s}^N, \notag
\end{align}
where $\lim_{N\rightarrow+\infty}\sup_{0\leq s\leq T}\mathbb{E}|\varepsilon_{16,s}^N|=0$. By Equations \eqref{equ 6.15}, \eqref{equ 6.17} and \eqref{equ 6.20},
\begin{align}\label{equ 6.21}
&G\left(\{V_{t,k}^N(f_k)\}_{k=0}^\mathcal{K}\right)-G(\{V_{0,k}^N(f_k)\}_{k=0}^\mathcal{K})-\int_0^t(\partial_s+\mathcal{L}_N)
G(\{V_{s,k}^N(f_k)\}_{k=0}^\mathcal{K})ds \\
&=G\left(\{V_{t,k}^N(f_k)\}_{k=0}^\mathcal{K}\right)-G(\{V_{0,k}^N(f_k)\}_{k=0}^\mathcal{K}) \notag\\
&\text{\quad}-\sum_{k=0}^{\mathcal{K}}\sum_{m=0}^{\mathcal{K}}\sum_{l=0}^{\mathcal{K}}\int_0^t\partial_kG(\{V_{s,k}^N(f_k)\}_{k=0}^\mathcal{K})
\left(V_{s,m}^N(\mathcal{P}_{s,m,l}^{k,1}f_k)+V_{s,l}^N(\mathcal{P}_{s,m,l}^{k,2}f_k)\right)ds \notag\\
&\text{\quad}-\frac{1}{2}\sum_{k_1=0}^\mathcal{K}\sum_{k_2=0}^\mathcal{K}\sum_{m=0}^\mathcal{K}\sum_{l=0}^\mathcal{K}\int_0^t\int_{\mathbb{T}^2}\partial^2_{k_1k_2}
G(\{V_{s,k}^N(f_k)\}_{k=0}^\mathcal{K})(b_{s,m,l}^{k_1}f_{k_1}(u,v))(b_{s,m,l}^{k_2}f_{k_2}(u,v))dudvds \notag\\
&\text{\quad}+\varepsilon_{17,t}^N, \notag
\end{align}
where $\lim_{N\rightarrow+\infty}\sup_{0\leq t\leq T}\mathbb{E}|\varepsilon_{17,t}^N|=0$. As we have shown in the proof Lemma \ref{lemma 6.2},
\[
\mathbb{E}\left((V_{t,k}^N(f))^2\right)=O(1).
\]
Hence, by Equation \eqref{equ 6.21},
\[
\left\{G\left(\{V_{t,k}^N(f_k)\}_{k=0}^\mathcal{K}\right)-G(\{V_{0,k}^N(f_k)\}_{k=0}^\mathcal{K})-\int_0^t(\partial_s+\mathcal{L}_N)G(\{V_{s,k}^N(f_k)\}_{k=0}^\mathcal{K})ds \right\}_{N\geq 1}
\]
are uniformly integrable for any $0\leq t\leq T$. Consequently, let $N\rightarrow+\infty$ in Equation \eqref{equ 6.21}, then
\begin{align*}
&\Bigg\{G\left(\{\widehat{V}_{t,k}(f_k)\}_{k=0}^\mathcal{K}\right)-G(\{\widehat{V}_{0,k}(f_k)\}_{k=0}^\mathcal{K}) \notag\\
&\text{\quad}-\sum_{k=0}^{\mathcal{K}}\sum_{m=0}^{\mathcal{K}}\sum_{l=0}^{\mathcal{K}}\int_0^t\partial_kG(\{\widehat{V}_{s,k}(f_k)\}_{k=0}^\mathcal{K})
\left(\widehat{V}_{s,m}(\mathcal{P}_{s,m,l}^{k,1}f_k)+\widehat{V}_{s,l}(\mathcal{P}_{s,m,l}^{k,2}f_k)\right)ds \notag\\
&\text{\quad}-\frac{1}{2}\sum_{k_1=0}^\mathcal{K}\sum_{k_2=0}^\mathcal{K}\sum_{m=0}^\mathcal{K}\sum_{l=0}^\mathcal{K}\int_0^t\int_{\mathbb{T}^2}\Big(
\partial^2_{k_1k_2}G(\{\widehat{V}_{s,k}(f_k)\}_{k=0}^\mathcal{K})\\
&\text{\quad}\times(b_{s,m,l}^{k_1}f_{k_1}(u,v))(b_{s,m,l}^{k_2}f_{k_2}(u,v))\Big)dudvds\Bigg\}_{0\leq t\leq T}
\end{align*}
is a martingale. In conclusion, $\{\widehat{V}_{t,k}: 0\leq t\leq T\}_{0\leq k\leq \mathcal{K}}$ is a solution to the martingale problem with respect to Equation \eqref{equ O-U} and the proof is complete.

\qed

\section{Applications} \label{section seven}
In this section, we apply our main results in exclusion processes and zero range processes.

\textbf{Application 1} \emph{Exclusion processes}. In the exclusion process, $\mathcal{K}=1$ and hence each position has at most one particle. For any $t\geq 0$, let
\[
\mu_t^N(du)=\frac{1}{N}\sum_{i=1}^N1_{\{\eta_t^N(i)=1\}}\delta_{i/N}(du)=\frac{1}{N}\sum_{i=1}^N\eta_t^N(i)\delta_{i/N}(du),
\]
then, under Assumption (A) for some $\psi\in\mathcal{A}_1$, Theorem \ref{theorem 2.1 mean field limit finite} shows that
\[
\lim_{N\rightarrow+\infty}\mu_t^N(f)=\int_{\mathbb{T}}f(u)\rho_t(u)du
\]
in $L^2$ for any $t\geq 0$ and $f\in C(\mathbb{T})$, where $\{\rho_t(u)\}_{t\geq 0, u\in \mathbb{T}}$ is the unique solution to
\begin{equation}\label{equ 6.1 EP mean field ode}
\begin{cases}
&\frac{d}{dt}\rho_t(u)=-\rho_t(u)\int_{\mathbb{T}}\phi_{1,0}(u,v)(1-\rho_t(v))dv\\
&\text{\quad\quad\quad\quad~}+(1-\rho_t(u))\int_{\mathbb{T}}\phi_{1,0}(v,u)\rho_t(v)dv \text{\quad for all~}u\in \mathbb{T},\\
&\rho_0=\psi.
\end{cases}
\end{equation}
Furthermore, let $V_t^N(du)=\frac{1}{\sqrt{N}}\sum_{i=1}^N\left(\eta_t^N(i)-\mathbb{E}\eta_t^N(i)\right)\delta_{i/N}(du)$, then Theorem \ref{theorem 2.3 fluctuation limit} shows that
\[
\lim_{N\rightarrow+\infty}V_t^N(f)=V_t(f)
\]
in distribution for any $t\geq 0$ and $f\in C^\infty(\mathbb{T})$, where $\{V_t\}_{t\geq 0}$ is the solution to
\begin{equation}\label{equ 6.3 EP fluctuation}
dV_t=\mathcal{P}_t^*dV_t+b_t^*d\mathcal{W}_t
\end{equation}
with $V_0(f)$ following the normal distribution $\mathbb{N}(0, \int_{\mathbb{T}}\psi(u)(1-\psi(u))f^2(u)du)$ for all $f\in C^\infty(\mathbb{T})$, where
\[
\mathcal{P}_tf(u)=\int_\mathbb{T}\phi_{1,0}(u,v)(1-\rho_t(v))(f(v)-f(u))dv-\int_{\mathbb{T}}\phi_{1,0}(u,v)\rho_t(v)(f(u)-f(v))dv
\]
and
\[
b_tf(u,v)=\sqrt{\phi_{1,0}(u,v)\rho_t(u)(1-\rho_t(v))}(f(v)-f(u))
\]
for all $f\in C^\infty(\mathbb{T})$ and $u,v\in \mathbb{T}$.

If $\phi_{1,0}$ is symmetric, i.e., $\phi_{1,0}(u,v)=\phi_{1,0}(v,u)$ for any $u,v\in \mathbb{T}$, then Equation \eqref{equ 6.1 EP mean field ode} reduces to the $C(\mathbb{T})$-valued linear ordinary differential equation
\begin{equation}\label{equ 6.2 SEP linear mean field equation}
\begin{cases}
&\frac{d}{dt}\rho_t(u)=\int_{\mathbb{T}}\phi_{1,0}(u,v)(\rho_t(v)-\rho_t(u))dv\text{\quad for all~}u\in \mathbb{T}, \\
&\rho_0=\psi
\end{cases}
\end{equation}
and Equation \eqref{equ 6.3 EP fluctuation} reduces to
\begin{equation}\label{equ 6.4 SEP fluctuation}
dV_t=\mathcal{P}^*dV_t+b_t^*d\mathcal{W}_t,
\end{equation}
where
\[
\mathcal{P}f(u)=\int_{\mathbb{T}}\phi_{1,0}(u,v)(f(v)-f(u))dv
\]
for all $f\in C^\infty(\mathbb{T})$ and $u\in \mathbb{T}$.

When $\phi_{1,0}$ is symmetric, the corresponding symmetric exclusion process is an example of the $N$-urn linear system introduced in \cite{Xue2023b}, the hydrodynamic limit and the fluctuation of which being driven by \eqref{equ 6.2 SEP linear mean field equation} and \eqref{equ 6.4 SEP fluctuation} respectively can also be deduced from Theorems 2.2 and 2.6 of \cite{Xue2023b}.

\quad

\textbf{Application 2} \emph{Zero range processes}. In the zero range process, $\mathcal{K}=+\infty$ and $\phi_{k,l}=\phi_k$ for all $k,l\geq 0$. Under Assumptions (A) and (B), according to the fact that the total number of particles is conserved, Theorem \ref{theorem 2.2 mean field limit infinite} shows that
\[
\lim_{N\rightarrow+\infty}\mu_{t,k}^N(f)=\int_{\mathbb{T}}\rho_{t,k}(u)f(u)du
\]
in $L^2$ for any $t\geq 0$, $k\geq 0$ and $f\in C(\mathbb{T})$, where $\{\rho_{t,k}(u): u\in \mathbb{T}, t\geq 0\}_{k\geq 0}$ is the solution to
\begin{equation}\label{equ zero range mean field limit}
\begin{cases}
&\frac{d}{dt}\rho_{t,k}(u)=-\rho_{t,k}(u)\int_{\mathbb{T}}\phi_k(u,v)dv-\rho_{t,k}(u)\int_{\mathbb{T}}\sum_{l=1}^{+\infty}\phi_l(v,u)\rho_{t,l}(v)dv\\
&\text{\quad\quad\quad\quad\quad}+\rho_{t,k-1}(u)\int_{\mathbb{T}}\sum_{l=1}^{+\infty}\phi_l(v,u)\rho_{t,l}(v)dv+\rho_{t, k+1}(u)\int_{\mathbb{T}}\phi_{k+1}(u,v)dv\\
&\text{\quad\quad\quad\quad\quad for all~}u\in \mathbb{T} \text{~and~} k\geq 1,\\
&\frac{d}{dt}\rho_{t,0}(u)=\rho_{t,1}(u)\int_{\mathbb{T}}\phi_1(u,v)dv-\rho_{t,0}(u)\int_{\mathbb{T}}\sum_{l=1}^{+\infty}\phi_l(v,u)\rho_{t,l}(v)dv\\
&\text{\quad\quad\quad\quad\quad for all~}u\in \mathbb{T},\\
&\rho_{0,k}(u)=e^{-\psi(u)}\frac{\psi(u)^k}{k!} \text{\quad for~}k\geq 0\text{~and~}u\in \mathbb{T},
\end{cases}
\end{equation}
and satisfies
\[
\limsup_{k\rightarrow+\infty}\frac{1}{k}\log \sup_{u\in \mathbb{T}}\left(\sum_{l=k}^{+\infty}l\rho_{t,l}(u)\right)=-\infty
\]
for all $t\geq 0$.

When $\phi_k=k\phi$ for any $k\geq 0$ and some $\phi\in C(\mathbb{T}^2)$, the zero range reduces to the generalized Ehrenfest model investigated in \cite{Xue2023}, where particles perform independent random walks. In this case, Equation \eqref{equ zero range mean field limit} reduces to
\begin{equation}\label{equ Ehrenfest mean field}
\begin{cases}
&\frac{d}{dt}\rho_{t,k}(u)=-k\rho_{t,k}(u)\int_{\mathbb{T}}\phi(u,v)dv-\rho_{t,k}(u)\int_{\mathbb{T}}\phi(v,u)\sum_{l=0}^{+\infty}l\rho_{t,l}(v)dv\\
&\text{\quad\quad\quad\quad\quad}+\rho_{t,k-1}(u)\int_{\mathbb{T}}\phi(v,u)\sum_{l=0}^{+\infty}l\rho_{t,l}(v)dv+(k+1)\rho_{t, k+1}(u)\int_{\mathbb{T}}\phi(u,v)dv\\
&\text{\quad\quad\quad\quad\quad for all~}u\in \mathbb{T} \text{~and~} k\geq 1,\\
&\frac{d}{dt}\rho_{t,0}(u)=\rho_{t,1}(u)\int_{\mathbb{T}}\phi(u,v)dv-\rho_{t,0}(u)\int_{\mathbb{T}}\phi(v,u)\sum_{l=0}^{+\infty}l\rho_{t,l}(v)dv\\
&\text{\quad\quad\quad\quad\quad for all~}u\in \mathbb{T},\\
&\rho_{0,k}(u)=e^{-\psi(u)}\frac{\psi(u)^k}{k!} \text{\quad for~}k\geq 0\text{~and~}u\in \mathbb{T}.
\end{cases}
\end{equation}
Let $\theta_t=\sum_{l=0}^{+\infty}l\rho_{t,l}$. Using Equation \eqref{equ Ehrenfest mean field}, we have
\[
\begin{cases}
&\frac{d}{dt}\theta_t(u)=-\theta_t(u)\int_{\mathbb{T}}\phi(u,v)dv+\int_{\mathbb{T}}\phi(v,u)\theta_t(v)dv, \\
&\theta_0(u)=\psi(u)
\end{cases}
\]
for all $u\in \mathbb{T}$. By Lemma \ref{lemma 5.1}, we have
\[
\lim_{M\rightarrow+\infty}\sup_{N\geq 1}\left(\frac{1}{N}\sum_{i=1}^N\sum_{k\geq M}k\mathbb{P}\left(\eta_t^N(i)=k\right)\right)=0.
\]
As a result, for any $f\in C(\mathbb{T})$,
\begin{align}\label{equ Ehrenfest hydrodynamic}
\lim_{N\rightarrow+\infty}\frac{1}{N}\sum_{i=1}^N\eta_t^N(i)f(i/N)
&=\lim_{N\rightarrow+\infty}\sum_{k=0}^{+\infty}k\mu_{t,k}^N(f)\\
&=\sum_{k=0}^{+\infty}k\int_{\mathbb{T}}\rho_{t,k}(u)f(u)du=\int_{\mathbb{T}}\theta_t(u)f(u)du \notag
\end{align}
in probability. Equation \eqref{equ Ehrenfest hydrodynamic} is first given by Theorem 2.3 of \cite{Xue2023}. Our above analysis gives another proof of this hydrodynamic limit conclusion.

\quad

\textbf{Acknowledgments.} The author is grateful to financial supports from the National Natural Science Foundation of China with grant number 12371142 and  the Fundamental Research Funds for the Central Universities with grant number 2022JBMC039.

{}
\end{document}